\documentclass[secthm]{elsart}
\usepackage{amsmath, amscd, amssymb, latexsym}
\newcommand{\bbbn}{\mathbb{N}}
\newcommand{\bbbs}{\mathbb{S}}
\newcommand{\bbbx}{\mathbb{X}}
\newcommand{\bbby}{\mathbb{Y}}
\newcommand{\bbbz}{\mathbb{Z}}
\newcommand{\leqct}{\leq_{\rm ct}}
\newcommand{\geqct}{\geq_{\rm ct}}
\newcommand{\eqct}{=_{\rm ct}}
\newcommand{\infct}{<_{\rm ct}}
\newcommand{\supct}{>_{\rm ct}}
\newcommand{\words}{{\bf 2}^*}
\newcommand{\dom}{\textit{dom}}
\newcommand{\prefix}{\textit{prefix}}
\newcommand{\lexico}{\textit{lexico}}
\newcommand{\weak}{\textit{wk}}
\newcommand{\strong}{\textit{st}}
\newcommand{\Reg}{\textit{Reg}}
\newcommand{\D}{{\mathcal D}}
\newcommand{\Dd}{{{\mathcal D}^{\text rev}}}
\newcommand{\XD}{{\bbbx\to\D}}
\newcommand{\WD}{{\words\to\D}}
\newcommand{\PR[1]}{{PR^{#1}}}
\newcommand{\Rec[1]}{{Rec^{#1}}}
\newcommand{\maxr[1]}{{Max^{#1}_{{\textit Rec}}}}
\newcommand{\maxpr[1]}{{Max^{#1}_{{\textit PR}}}}
\newcommand{\minr[1]}{{Min^{#1}_{{\textit Rec}}}}
\newcommand{\minpr[1]}{{Min^{#1}_{{\textit PR}}}}
\newcommand{\kmax[1]}{{K^{#1}_{{\textit max}}}}
\newcommand{\kmin[1]}{{K^{#1}_{{\textit min}}}}
\begin{document}
\begin{frontmatter}
\title{Kolmogorov complexities $\kmax[]$, $\kmin[]$\\
on computable partially ordered sets}
\author{Marie Ferbus-Zanda}
\address{LIAFA, Universit\'e Paris 7 \& CNRS, France\\
{\tt ferbus@logique.jussieu.fr}}
\author{Serge Grigorieff}
\address{LIAFA, Universit\'e Paris 7 \& CNRS, France\\
{\tt seg@liafa.jussieu.fr}}
\tableofcontents
\begin{abstract}
We introduce a machine free mathematical framework to get a natural
formalization of some general notions of infinite computation
in the context of Kolmogorov complexity.
Namely, the classes $\maxpr[\XD]$ and $\maxr[\XD]$
of functions $\XD$ which are pointwise maximum of
partial or total computable sequences of functions
where $\D=(D,<)$ is some computable partially ordered set.
The enumeration theorem and the invariance theorem always hold
for $\maxpr[\XD]$, leading to a variant $\kmax[\D]$ of
Kolmogorov complexity.
We characterize the orders $\D$ such that the enumeration theorem
(resp. the invariance theorem) also holds for $\maxr[\XD]$.
It turns out that $\maxr[\XD]$ may satisfy the invariance theorem
but not the enumeration theorem.
Also, when $\maxr[\XD]$ satisfies the invariance theorem
then the Kolmogorov complexities associated to $\maxr[\XD]$ and
$\maxpr[\XD]$ are equal (up to a constant).
\\
Letting $\kmin[\D]=\kmax[\Dd]$, where $\Dd$ is the reverse order,
we prove that
either $\kmin[\D]\eqct\kmax[\D]\eqct K^D$
($\eqct$ is equality up to a constant)
or $\kmin[\D],\kmax[\D]$ are $\leqct$ incomparable and
$\infct K^D$ and $\supct K^{\emptyset',D}$.
We characterize the orders leading to each case.
We also show that $\kmin[\D],\kmax[\D]$ cannot be both
much smaller than $K^D$ at any point.
\\
These results are proved in a more general setting with
two orders on $D$, one extending the other.
\end{abstract}
\end{frontmatter}
%
\section{Introduction}
\label{s:intro}
%
\subsection{Non halting programs for which the current output
is eventually the wanted object, but one does not know when...}
\label{ss:eventually}
%
In this paper, we consider a particular kind of
description methods
in order to define variants of Kolmogorov complexity.
Let's start with two paradigmatic examples.
Given $n\in\bbbn$ and $u\in\Sigma^*$
(where $\Sigma$ be some finite alphabet),
how do we get
\\\indent- the value $BB(n)$ of the busy beaver function
           $BB:\bbbn\to\bbbn$,
\\\indent- the value $K_{\Sigma^*}(u)$ of Kolmogorov complexity
           $K_{\Sigma^*}:\Sigma^*\to\bbbn$ ?
\\
The definitions of $BB(n)$ and $K_{\Sigma^*}(u)$ lead to the
following mechanisms.
\\\indent- run all Turing machines with $\leq n$ states
           and all programs with length $\leq n$,
\\\indent- for each $t$, consider those machines and programs
           halting in $\leq t$ steps,
\\\indent- look at the maximum number of cells
           visited by these machines,
\\\indent- look at the minimum length of these programs.
\\
In this way, one gets two computable functions
$bb:\bbbn\times\bbbn\to\bbbn$ and $k:\Sigma^*\times\bbbn\to\bbbn$,
with one more integer argument (for time steps)
such that,
for every fixed $n\in\bbbn$ and $u\in\Sigma^*$, the maps
$t\mapsto bb(n,t)$ and $t\mapsto k(u,t)$ are respectively
monotone increasing and decreasing
and are both eventually constant with respective values
$BB(n)$ and $K_{\Sigma^*}(u)$.
Since neither $BB$ nor $K_{\Sigma^*}$ is computable, there is no
computable functions of $n$ or $u$ which bound the moment
these maps become constant.
\medskip\\
These examples lead us to introduce the following notion
of description methods for objects of a partially ordered
set $\D$ with a computable structure
(cf. Definition \ref{def:computableposet}, \ref{def:PRinfini}).
\\
{\em A computable approximation from below (resp. from above)
of objects of $\D$ is a program for a computable function
$f:\bbbx\times\bbbn\to\D$ (where $\bbbx$ is some reasonable set
such as $\bbbn$ or $\words$, cf. \S\ref{ss:notation} Notations)
such that,
for every fixed $x\in\bbbx$, the map $t\mapsto f(x,t)$ is
monotone increasing (resp. decreasing) and eventually constant.
Nothing is assumed about the moment $t\mapsto f(x,t)$ becomes
constant:
there may be no computable function of $x$ majorizing it.
\\
The associated decompressor
--- or description method ---
is the function $F:\bbbx\to\D$ such that
$F(x)$ is the limit value of $f(x,t)$ when $t\to+\infty$,
i.e. the maximum (resp. minimum) value of the finite set
$\{f(x,t):t\in\bbbn\}$.}
We shall call such functions $F$ computably approximable
from below (resp. from above).
\medskip\\
Consider $\Sigma^*$ with the prefix ordering.
The context of non halting (hence infinite) computations,
cf. Chaitin, 1975 \cite{chaitin75},
and Solovay, 1977 \cite{solovay77},
leads to functions $F:\words\to\Sigma^*$ which are computably
approximable from below.
In fact, if the output alphabet is $\Sigma$,
the current output $f(x,t)$ at time $t$ is a function
$f:\words\times\bbbn\to\Sigma^*$ which is a computable
approximation from below for words in $\Sigma^*$ such that $F(x)$
is the $\max$ of the $f(x,t)$'s .
\medskip\\
Observe that,
in case the ordered set $\D$ is noetherian (resp. well-founded),
the notion of approximation from below (resp. from above) of
objects of $\D$ reduces to that of computable function
$f:\bbbx\times\bbbn\to\D$ which is monotone increasing
(resp. decreasing) with respect to its second argument.
This is indeed the case with the approximation from above of
the values of $K_{\Sigma^*}$ since $(\bbbn,<)$ is well-founded.
Cf. also \S\ref{sss:noether}.
\medskip\\
Other examples are developped in \S\ref{ss:examples}.
In particular, there is one involving quotients of regular languages
by a fixed computably enumerable language.
%
%
\subsection{Functions approximable from below (resp. from above)
as decompressors for variants of Kolmogorov complexity}
\label{ss:functions}
%
The above mentioned context of non halting computations
has recently led to interesting variants
$K^\infty_{\Sigma^*}:\Sigma^*\to\bbbn$,
$K^\infty_\bbbn:\bbbn\to\bbbn$ of Kolmogorov complexity
introduced (in their prefix-complexity version $H^\infty$)
in Becher \& Chaitin \cite{becherchaitindaicz}.
\\
This last Kolmogorov complexity $K^\infty_\bbbn$ has also proved to
be equal to the Kolmogorov complexity $K_{\textit card}$ introduced
in Ferbus \& Grigorieff,
2002 \cite{ferbusgrigoClermont,ferbusgrigoKandAbstraction1}
where we compare some natural set theoretical semantics of integers,
namely
Church iterators of functions,
cardinals of computably enumerable sets,
indexes of computably enumerable equivalence relations.
Comparison of these semantics is done via associated Kolmogorov
complexities which somehow constitute measures of their
``abstraction degree" and are defined in terms of infinite
or/and oracular computations.
\medskip\\
The cornerstone of Kolmogorov complexity, namely the invariance
theorem, really deals with partial computable functions,
not Turing machines.
In fact, Turing machines do not constitute such an abstract
structured mathematical framework as partial computable functions do.
Going to this last framework opens new natural considerations
which would not be simply viewed with Turing machines.
\\
In this paper we abstract from non halting computations on
Turing machines and develop a general machine-free
mathematical framework using a partially ordered set $\D$.
Namely, letting $\bbbx$ be a basic space
(cf. \S\ref{ss:notation} Notations),
we introduce the classes of functions $F:\bbbx\to\D$
$$
\maxpr[\XD]\ \ ,\ \ \minpr[\XD]
$$
which are partial computably approximable from below
(resp. from above).
Which means that the $f:\words\times\bbbn\to\D$ such that
$F(x)$ is the $\max$ or $\min$ of the $f(x,t)$'s
is partial computable rather than computable.
\\
Of course, the ${\textit Min}$ classes are the ${\textit Max}$ classes
associated to the reverse order.
\medskip\\
We also introduce the subclasses of functions
$$
\maxr[\XD]\ \ ,\ \ \minr[\XD]
$$
which are computably approximable from below (resp. from above).
It happens that the $\maxr[\XD]$ class is closely related with
the class based on non halting Turing machines computations
with outputs in $\D$ (modulo adequate coding of $\D$).
\medskip\\
As for the above examples, the busy beaver function
$BB:\bbbn\to\bbbn$ is in $\maxr[\D]$and
Kolmogorov complexity $K_{\Sigma^*}:\Sigma^*\to\bbbn$ and
its prefix-free variant $H_{\Sigma^*}:\Sigma^*\to\bbbn$
are in $\minr[\D]$ with $\D=(\bbbn,<)$.
\medskip\\
These classes lead to new variants of Kolmogorov complexity
which would just be ignored when considering Turing machines.
%
\subsection{Main theorems}
\label{ss:main}
%
The development of Kolmogorov complexities $\kmax[\D]$, $\kmin[\D]$
associated to the classes $\maxpr[\WD]$ and $\minpr[\WD]$
is straightforward (cf. \S\ref{s:Kmax}).
The main results of the paper deal with the comparison of
$\kmax[\D]$, $\kmin[\D]$ with the classical Kolmogorov complexity
$K^D$ and its relativized version $K^{\emptyset',D}$ to oracle
$\emptyset'$.
In \S\ref{s:hierarchy}, we prove three theorems which give
the main comparison relations
(relative to the ``up to a constant" order $\leqct$,
cf. \S\ref{ss:notation} Notations)
between these complexities.
\medskip\\
The first theorem (Thm.\ref{thm:hierarchy1}) is valid
whatever be the partial order on $D$.
It states that $K^{\emptyset',D}\infct\inf(\kmax[\D],\kmin[\D])$
and that $\kmax[\D]$, $\kmin[\D]$, though obviously $\leqct K^D$,
cannot be simultaneously much smaller than $K^D$ since
$$
K^D\leqct(\kmax[\D]+\log(\kmax[\D]))+(\kmin[\D]+\log(\kmin[\D]))
$$
The second theorem (Thm.\ref{thm:hierarchy2}) proves that
either $\kmax[\D]\eqct\kmin[\D]\eqct K^D$
or     $\kmax[\D]$,$\kmin[\D]$ are $\leqct$ incomparable and
both are $\infct$ to $K^D$.
This dichotomy is also characterized by a simple property
on the order.
\medskip\\
The third theorem (Thm.\ref{thm:hierarchy3}) considers two partial
orders $<_\weak$ and $<_\strong$ on $D$, the second extending
the first.
We give conditions $(*)$ and $(**)$ on the orders such that
\\-
$(*)$ insures that $\kmax[\D_\strong]\eqct\kmax[\D_\weak]$
and $\kmin[\D_\strong]\eqct\kmin[\D_\weak]$,
\medskip\\-
$(**)$ insures that $\kmax[\D_\strong]\infct\kmax[\D_\weak]$
and $\kmax[\D_\strong]\infct\kmax[\D_\weak]$
and neither $\kmax[\D_\strong]$ nor $\kmin[\D_\strong]$
is $\leqct\min(\kmax[\D_\weak],\kmin[\D_\weak])$.
\medskip\\
These conditions are almost complementary: $(**)$ is an effective
version of the negation of $(*)$.
\\
An interesting case of this theorem is obtained when
$\Sigma=\{1,...,k\}$ with the obvious order and
$<_\weak$, $<_\strong$ are the prefix and the lexicographic
orders on $\Sigma^*$ (the last one being isomorphic to the order
on $k$-adic rational reals in $[0,1]$).
%
%
\subsection{The $\maxr[]$ and $\minr[]$ classes}
\label{ss:maxrminr}
%
In \S\ref{ss:PRversusMaxPR} and \ref{ss:syntax} we come back to the
four classes
$\maxpr[\WD]$, $\minpr[\WD]$ and $\maxr[\WD]$, $\minr[\WD]$.
We compare them to that of partial computable functions $\XD$
and look at the syntactical complexity of their domains and graphs.
\\
In \S\ref{ss:MaxInterMin} we compute $\maxpr[\XD]\cap\minpr[\XD]$
under simple conditions about the partial order on $D$.
\medskip\\
In \S\ref{s:KwithMaxRec}, we consider the possible development
of Kolmogorov complexities based on the classes
$\maxr[\WD]$ and $\minr[\WD]$.
This leads to look at the two following problems:
\\\indent- the existence of an enumeration,
\\\indent- the invariance theorem.
\\
For each problem, we characterize the orders $\D$ for which
there is a positive answer
(cf. \S\ref{ss:enumMaxRec}, \ref{ss:KmaxRec}).
\\
It turns out (cf. Thm.\ref{thm:KmaxRec}) that when the invariance
theorem holds for $\maxr[\WD]$ then every function in $\maxpr[\WD]$
has an extension (not necessarily total) in $\maxr[\WD]$.
This insures that the Kolmogorov complexities associated to
$\maxr[\WD]$ and $\maxpr[\WD]$ coincide.
In particular, $K^\infty_{\Sigma^*}$ is the complexity associated to
$\maxr[\WD]$ and $\maxpr[\WD]$ when $\D$ is $\Sigma^*$ with the
prefix order.
\\
Surprisingly, there are orders such that the invariance
theorem holds for $\maxr[\WD]$ whereas the enumeration theorem fails
(compare Thm.\ref{thm:enumMaxr} and Thm.\ref{thm:KmaxRec}).
%
%
\subsection{Notations}
\label{ss:notation}
%
1. Equality, inequality and strict inequality up to a constant
between total functions $\bbbs\to\bbbn$ are denoted as follows:
\medskip\\
$\begin{array}{rcccl}
f\ \leqct\ g &\Leftrightarrow &
\exists c\ \forall {\tt s}\ f({\tt s})\leq g({\tt s})+c&&\medskip\\
f\ \eqct\ g &\Leftrightarrow & f\leqct g\ \wedge\ g\leqct f
&\Leftrightarrow &
\exists c\ \forall {\tt s}\ |f({\tt s})-g({\tt s})|\leq c\medskip\\
f\ \infct\ g &\Leftrightarrow &
f \leqct g\ \wedge\ \neg(g \leqct f)
&\Leftrightarrow &f \leqct g\
\wedge\ \forall c\ \exists {\tt s}\ g({\tt s})>f({\tt s})+c
\end{array}$
\medskip\\
2. {\em [Basic spaces]\ }
$\words$ denotes the set of binary words.
We call basic spaces the products of non empty finite families
of spaces of the form $\bbbn$ or $\bbbz$ or $A^*$ where $A$ is
some finite alphabet.
Basic spaces are denoted by $\bbbs,\bbbx,\bbby,...$
\medskip\\
3.  {\em [Partial recursive (or computable) functions]\ }
$\PR[\bbbx\to\bbby]$ (resp. $Rec[\bbbx\to\bbby]$) denotes
the family of partial (resp. total) computable functions from
$\bbbx$ to $\bbby$.
%
%
\section{The $Max$ and $Min$ classes of functions}
\label{s:MaxMin}
%
%
\subsection{Infinite computations and monotone machines}
\label{ss:MachinePRinfini}
%
Recall that a Turing machine is monotone if its current
output may only increase with respect to the prefix
order on words: no overwriting is allowed.
This is indeed Turing's original assumption \cite{turing36},
insuring that, in the limit of time, the output of a non halting
computation always converges,
either to a finite or to an infinite sequence.
This concept was also considered by Levin
\cite{levin73} and Schnorr \cite{schnorr73,schnorr77},
see \cite{livitanyi} p.276.
Such infinite computations with possibly infinite outputs
can be used to obtain highly random reals,
cf. Becher \& Chaitin \cite{becherchaitindaicz}
and Becher \& Grigorieff \cite{bechergrigo}.
\\
{\em In this paper, when considering infinite computations,
we retain the sole limit outputs that are finite.}
\medskip\\
The following easy proposition links infinite computations, as
considered for the definition of $K^\infty$ and its prefix version
$H^\infty$ introduced in \cite{bechernies},
with the general approach which is the subject of this paper.
\begin{prop}\label{p:machineKI}
Let $F:\words\to\Sigma^*$ where $\Sigma$ is some non empty finite
alphabet.
The following conditions are equivalent:
\begin{enumerate}
\item[i.]
$F$ can be computed via possibly infinite computations on some
monotone Turing machine with output alphabet $\Sigma$, according to
the following convention:
$F({\tt s})$ is defined if and only if the output remains constant
after some step.
\item[ii.]
There exists a total computable function
$f:\words\times\bbbn\to\Sigma^*$
such that
\\ - $f({\tt s},t)$ is monotone increasing in $t$ with respect to the
prefix order on $\Sigma^*$,
\\ - ${\tt s}\in \dom(F)$ if and only if
$\{f({\tt s},t):t\in\bbbn\}$ is finite and non empty,
\\ - $F({\tt s})$ is the maximum value of
$\{f({\tt s},t):t\in\bbbn\}$.
\item[iii.]
Let $\lambda$ denote the empty word.
Idem as ii, with $f$ such that
$$
f({\tt s},0)=\lambda\ \ \ ,\ \ \
f({\tt s},t+1)\in\{f({\tt s},t)\}\cup\{f({\tt s},t)\sigma:\sigma\in\Sigma\}
$$
\end{enumerate}
\end{prop}
\begin{pf}
$iii\Rightarrow ii$ is trivial;
$i\Leftrightarrow iii$ :
let $f({\tt s},t)$ be the current output at time $t$ when the input is
${\tt s}$.
As for $ii\Rightarrow iii$, let
$\widetilde{f}({\tt s},0)=\lambda$
and
$\widetilde{f}({\tt s},t+1)$ be the prefix of $f({\tt s},t+1)$
with length $\min(|\widetilde{f}({\tt s},t)|+1,|f(\tt s,t+1)|)$.
Then
$\{\widetilde{f}({\tt s},t):t\in\bbbn\}$ and
$\{f({\tt s},t):t\in\bbbn\}$ are simultaneously finite or infinite
and, when finite, their maximum elements are equal.
\qed\end{pf}
%
\subsection{Mathematical modelization:
the $Max$ and $Min$ classes}
\label{ss:Max}
%
Proposition \ref{p:machineKI} and the argumentation in
\S\ref{ss:eventually}--\ref{ss:functions}
invite to a mathematical, machine-free modelization of
the notion of function defined by infinite computations.
Namely that of function obtained as pointwise maximum of
a computable sequence of total computable functions.
{\em A construction which makes sense for maps from a basic
set $\bbbx$ into any computable partially ordered set $\D=(D,<)$},
and leads to the class $\maxr[\XD]$.
\\
It is also quite natural
-- {\em in fact, it is even much more natural
from a mathematical point of view} --
to consider the version of the above modelization using
{\em partial computable functions} instead of total computable ones.
This leads to the class $\maxpr[\XD]$.
\\
Natural and interesting important examples (cf. \S\ref{ss:examples})
are obtained when $\D$ is among the following (obviously computable)
partially ordered sets:
$$
(\bbbn,<)\ \ ,\ \ (\bbbz,<)\ \ ,\ \
(\Sigma^*,<_{\textit prefix})\ \ ,\ \ (\Sigma^*,<_{\textit lexico})
$$
and the reverse orders obtained by replacing $<$ by $>$,
where $<_{\textit lexico}$ on $\Sigma^*$ depends on a total or
partial order on the alphabet $\Sigma$.
\begin{defn}[The $\max^\D$ and $\min^\D$ operators]
Let $\bbbx$ be some basic set and
${\mathcal D}=(D,<)$ be some partially ordered set.
Let $f:\bbbx\times\bbbn\to D$ be monotone increasing in its second
argument on its domain.
We define ${\max}^\D f:\bbbx\to D$ (resp. ${\min}^\D f:\bbbx\to D$)
as the function
\begin{enumerate}
\item[i.]
defined on the ${\tt x}$'s in $\bbbx$ for which the map
$t\mapsto f({\tt x},t)$ has finite non empty range,
\item[ii.]
and such that
$({\max}^\D f)({\tt x})$ (resp. $({\min}^\D f)({\tt x})$)
is the maximum (resp. minimum) element of
$\{f({\tt x},t):t\in\bbbn\}$.
\end{enumerate}
\end{defn}
\begin{defn}\label{def:computableposet}$\\ $
1. A computable partially ordered set ${\mathcal D}$ is a triple
$(D,<,\rho)$ such that $\rho:\bbbn\to D$ is a bijective total map
(in particular, $D$ is infinite countable)
and $<$ is a partial order on $D$ such that
$\{(m,n):\rho(m)<\rho(n)\}$ is computable.
\\
2. Let $\bbbx$ be a basic space.
A function $F:\bbbx\to D$ is partial (resp. total) computable
if so is $\rho^{-1}\circ F:\bbbx\to\bbbn$.
\\
A set $Z\subseteq\bbbx\times D^k$ is computable if so is
$(Id_\bbbx,\rho,...,\rho)^{-1}(Z)$ as a subset of
$\bbbx\times\bbbn^k$, where $Id_\bbbx$ is the identity function
on $\bbbx$.
\end{defn}
Of course, we shall omit any reference to $\rho$ when $\D$ is
$\bbbn$ or $\bbbz$ with the natural order,
or $\Sigma^*$ with the prefix or the lexicographic order
(with respect to some partial or total order of the elements
of $\Sigma$).
\begin{defn}[$Max$ and $Min$ classes]
\label{def:PRinfini}
Let $\bbbx$ be a basic space and $\D=(D,<,\rho)$ be a computable
partially ordered set. We let
\begin{eqnarray*}
\maxr[\XD]&=&\{{\max}^\D f : f:\bbbx\times\bbbn\to D
\mbox{ is total computable}\}
\\
\maxpr[\XD]&=&\{{\max}^\D f : f:\bbbx\times\bbbn\to D
\mbox{ is partial computable}\}
\end{eqnarray*}
We respectively denote by $\minpr[\XD]$ and $\minr[\XD]$ the analog
classes defined with the ${\min}^\D$ operator, i.e. the classes
$\maxpr[\bbbx\to\Dd]$ and $\maxpr[\bbbx\to\Dd]$ where $\Dd=(D,>)$.
\end{defn}
Proposition \ref{p:machineKI} can be rephrased in terms of the
prefix ordering on $\Sigma^*$.
\begin{prop}\label{p:exitMach}
If $\Sigma$ is a finite alphabet then
$\maxr[\words\to(\Sigma^*,<_\prefix)]$
is the class of functions computed via possibly infinite
computations on monotone Turing machines
(cf. Proposition \ref{p:machineKI} i)
with $\Sigma$ as output alphabet.
\end{prop}
%
%
\subsection{Domains of functions in the $Max/Min$ classes}
\label{ss:domains}
%
We denote by $\Sigma^0_1\wedge\Pi^0_1$ the family of conjunctions
of $\Sigma^0_1$ and $\Pi^0_1$ formulas.
Let $\bbbx$ be a basic set and $\D$ be a computable ordered set.
The arithmetical hierarchy on $\bbbn$ induces a hierarchy on $D$ and
$\bbbx\times D$ :
a relation $R\subseteq\bbbx\times D$ is
$\Sigma^0_n$ or $\Pi^0_n$ or $\Sigma^0_n\wedge\Pi^0_n$
if so is $(Id_\bbbx,\rho)^{-1}(R)\subseteq\bbbx\times\bbbn$.
\begin{prop}\label{p:syntax}
Let $\bbbx$ be a basic set and $\D$ be a computable ordered set.
\\
Every partial function in $\maxpr[\XD]$ or in $\minpr[\XD]$ has
$\Sigma^0_1\wedge\Pi^0_1$ graph and $\Sigma^0_2$ domain.
\end{prop}
\begin{pf}
Let $f:\bbbx\times\bbbn\to D$ be partial computable,
monotone increasing in its second argument on its domain. Then
\begin{eqnarray*}
({\max}^\D f)({\tt x})=z&\Leftrightarrow&
\exists t\ (f({\tt x},t)
\mbox{ is defined }\wedge\ f({\tt x},t)=z)
\\ && \wedge\ \forall t\
(f({\tt x},t)\mbox{ is defined }
               \Rightarrow\ f({\tt x},t)\leq z)
\medskip\\
{\tt x}\in \dom({\max}^\D f)&\Leftrightarrow&\exists z\ F({\tt x})=z
\end{eqnarray*}
Idem with $\minpr[\XD]$.
\qed\end{pf}
%
%
\subsection{Examples of functions in the $Max$ and $Min$ classes}
\label{ss:examples}
%
The classes $\maxr[\XD],\minr[\XD]$ contain
many fundamental non computable functions.
To see that some functions are not in such classes,
we shall use Theorem \ref{thm:MaxInterMin} below
(the proof of which does not depend on any result of this \S).
%
\subsubsection{Kolmogorov and Chaitin-Levin program-size
complexities}
\label{sss:exK}
%
\begin{prop}\label{p:exK}
Let $\D$ be $(\bbbn,<)$.
Kolmogorov and Chaitin-Levin program-size complexities
$K_\bbbn,H_\bbbn:\bbbn\to\bbbn$
(resp.$K_{\Sigma^*},H_{\Sigma^*}:\Sigma^*\to\bbbn$)
are in $\minr[\bbbn\to\D]\setminus\maxpr[\bbbn\to\D]$
(resp. in
$\minr[\Sigma^*\to\D]\setminus\maxpr[\Sigma^*\to\D]$).
\end{prop}
\begin{pf}
That $K,H$ belong to $\minr[\bbbn\to\D]$ is a mere
reformulation of the well-known fact that they are
computably approximable from above, i.e. they are limits of
decreasing computable sequences of total computable functions.
That these {\em total} functions are not in $\maxpr[\bbbn\to\D]$ is
an obvious application of Theorem \ref{thm:MaxInterMin} below.
\qed\end{pf}
%
%
\subsubsection{Busy beaver}\label{sss:exBB}
%
\begin{prop}\label{p:exBB}
Let $\D$ be $(\bbbn,<)$.
Let $BB:\bbbn\to\bbbn$ be the busy beaver function, i.e.
$BB(n)$ is the maximum number of cells visited by the input head
of a Turing machine with $n+1$ states which halts with no input.
\\
Then $BB\in\maxr[\bbbn\to\D]\setminus\minpr[\bbbn\to\D]$.
\end{prop}
\begin{pf}
Observe that $BB=\max bb$ where $bb$ is the total computable
function such that  $bb(n,t)$ is the maximum among $0$ and
the numbers of cells visited by Turing machines with $n+1$
states which halt in at most $t$ steps.
\\
An obvious application of Theorem \ref{thm:MaxInterMin} below
shows that $BB$ is not in $\minpr[\bbbn\to\D]$.
\qed\end{pf}
\begin{rem}
Variants of the busy beaver function can be very
naturally defined with ranges over various types of data structures.
For instance, finite graphs relative to the inclusion or embedding
ordering.
\end{rem}
%
\subsubsection{Cardinality of finite computably enumerable sets}
\label{sss:excardRE}
%
The following example is completely investigated in
\cite{ferbusgrigoClermont,ferbusgrigoKandAbstraction1}.
\begin{prop}\label{p:excardRE}
Let $\D$ be $(\bbbn,<)$.
Let $\textit{cardRE}:\bbbn\to\bbbn$ be such that
\begin{eqnarray*}
\textit{cardRE}(n)&=&\left\{
\begin{array}{ll}
card(W_n)&\mbox{if $W_n$ is finite}
\\\mbox{undefined}&\mbox{otherwise}
\end{array}\right.
\end{eqnarray*}
where $card(W_n)$ is the number of elements of the
computably enumerable set $W_n$ with code $n$.
\\
Then
$\textit{cardRE}\in\maxr[\bbbn\to\D]\setminus\minpr[\bbbn\to\D]$.
\end{prop}
\begin{pf}
Observe that $\textit{cardRE}=\max h$ where $h(n,t)$ is total
computable and counts the number of elements of $W_n$ obtained after
$t$ computation steps.
\\
The domain of the partial function $\textit{cardRE}$ is known to be
$\Sigma^0_2$ complete, hence not $\Sigma^0_1\wedge\Pi^0_1$.
Applying Theorem \ref{thm:MaxInterMin} below, we see that
$\textit{cardRE}$ cannot be in $\minpr[\bbbn\to\D]$.
\qed\end{pf}
%
%
\subsubsection{Interacting finite sets with a fixed computably
enumerable set}
\label{sss:exFinRE}
%
If $X,Y\subseteq\bbbn$, let's denote $X-Y$ and $X\setminus Y$
the sets
$$
X-Y=\{x-y:x\in X\wedge y\in Y\wedge x\geq y\}\ ,\
X\setminus Y=\{z:z\in X\wedge z\notin Y\}
$$
\begin{prop}\label{p:interacting}
Let $\D$ be be the family $P_{<\omega}(\bbbn)$ of finite subsets
of $\bbbn$, ordered by set inclusion.
If $A\subseteq\bbbn$ is a fixed computably enumerable set which is
non computable then
\begin{enumerate}
\item[i.]
the maps $X\mapsto X\cap A$ and $X\mapsto X-A$
are in $\maxr[\D\to\D]\setminus\minpr[\D\to\D]$.
\item[ii.]
the map $X\mapsto X\setminus A$ is in
$\minr[\D\to\D]\setminus\maxpr[\D\to\D]$.
\end{enumerate}
\end{prop}
\begin{pf}
Let $A=\varphi(\bbbn)$ where $\varphi:\bbbn\to\bbbn$ is total
computable.
Define total computable maps $f,g,h:\D\times\bbbn\to D$ such that
\medskip\\\medskip\centerline{$\begin{array}{rclcrcl}
f(X,t)&=&X\cap\varphi(\{0,...,t\})&
\hspace{5mm}g(X,t)&=&X-\varphi(\{0,...,t\})
\\
h(X,t)&=&X\setminus\varphi(\{0,...,t\})
\end{array}$}
It is easy to see that $X\cap A=({\max}^\D f)(X)$ and
$X- A=({\max}^\D g)(X)$ and $X\setminus A=({\min}^\D h)(X)$ .
\qed\end{pf}
%
%
\subsubsection{Quotients of regular languages by a fixed computably
enumerable language}
\label{sss:exQuotient}
%
We now come to a very different example.
\medskip\\
The family $\Reg$ of regular languages over alphabet $\Sigma$
can be defined by regular expressions which are words in the
alphabet $\widetilde{\Sigma}$ obtained by enriching $\Sigma$ with
symbols $+,*,\cdot,(,)$.
\\
Let $\zeta:\widetilde{\Sigma}^*\to\Reg\,$ be the surjective map
such that,
if $u$ is a regular expression then $\zeta(u)$ is the associated
regular language,
else $\zeta(u)=\emptyset$.
\\
Since equality of regular languages is decidable, there exists a
computable map $\eta:\bbbn\to\widetilde{\Sigma}^*$ such that
$\rho=\zeta\circ\eta:\bbbn\to\Reg\,$ is bijective.
\\
Using decidability of inclusion of regular languages, we see that
$(\Reg,\subseteq,\rho)$ is a computable partially ordered set
in the sense of Definition \ref{def:computableposet}.
\medskip\\
It is known that, if $L$ is a regular language and
$M\subseteq\Sigma^*$ is any language (even non computable) then
$$
M^{-1}L=\{u\in\Sigma^*:\exists v\in M\ vu\in L\}
$$
is always regular and $M^{-1}L=M'^{-1}L$ for some finite subset
$M'\subseteq M$.
Recall the core of the easy proof:
if $L$ is the set of words leading from state $q_0$ to a final state
of automaton $\mathcal A$
and if the words in $M$ lead from state $q_0$ to the states in $X$,
then $M^{-1}L$ is the set of words leading from a state in $X$
to a final state.
\begin{prop}\label{p:exQuotient}
Let $M\subseteq\Sigma^*$ be a fixed computably
enumerable language which is non computable.
Let $F_M:\Reg\to\Reg\,$ be such that $F_M(L)=M^{-1}L$.
Then $F_M$ is in
$\maxr[\Reg\to\Reg]\setminus\minpr[\Reg\to\Reg]$.
\end{prop}
\begin{pf}
Let $M=\varphi(\bbbn)$ where $\varphi:\bbbn\to\Sigma^*$
is a total computable function.
Observe that $F_M=\max^{\Reg}f_M$ where $f_M:\Reg\times\bbbn\to\Reg$
is such that
$$
f_M(L,t)=(\varphi(\{0,...,t\})^{-1}L
$$
Observe that $M$ is computable with oracle $F$ since
$u\in M$ if and only if $M^{-1}\{u\}=\{\lambda\}$.
Since $M$ is not computable, $F$ cannot be computable.
Using Theorem \ref{thm:MaxInterMin} point 1
(and the fact that $F$ is total), we see that $F$ is not in
$\minpr[\Reg\to\Reg]$.
\qed\end{pf}
Using the above surjection $\zeta:\widetilde{\Sigma}^*\to\Reg\,$,
one can reformulate the above result in terms of a partial
computable preordering on words quite different of the usual ones.
This necessitates a straightforward extension to preorderings of the
material about the $Max$ and $Min$ classes.
\\
Let $\mu:\Reg\to\Delta^*$ be the map which associates to a regular
language $L$ the regular expression
(obtained via some fixed algorithm)
describing its minimal automaton.
Observe that $\zeta$ is a retraction of the injective map $\mu$,
i.e. $\zeta\circ\mu$ is the identity map on $\Reg\,$.
\begin{prop}\label{p:exQuotient2}
Let $\D$ be $\widetilde{\Sigma}^*$ with the following computable
preordering:
$$
u\preceq v\ \Leftrightarrow\ \zeta(u)\subseteq\zeta(v)
$$
Let $M\subseteq\widetilde{\Sigma}^*$ be a fixed computably enumerable
language which is non recucomputablersive.
Then the map $u\mapsto\mu(M^{-1}\zeta(u))$
(which maps a regular expression for $L$ to one for $M^{-1}L$)
is in $\maxr[\D\to\D]\setminus\minpr[\D\to\D]$.
\end{prop}
\begin{pf}
Let $F,f$ be as in the proof of Proposition \ref{p:exQuotient}.
Since $\zeta\circ\mu=Id_\Reg$, we see that
$\widetilde{F}=\mu\circ F\circ\zeta$ makes the following
diagram commute:
\[ \begin{CD}
\Reg                  @> F >>              \Reg\\
@A\zeta               AA @AA               \zeta A \\
\widetilde{\Sigma}^*  @> \widetilde{F} >>  \widetilde{\Sigma}^*
\end{CD} \]
which allows to transfer the results of
Proposition \ref{p:exQuotient}.
\qed\end{pf}
%
%
\subsubsection{Noetherian or well-founded orderings}
\label{sss:noether}
%
Suppose $\D$ is Noetherian (resp. well-founded) and let
$f:\bbbx\times\bbbn\to D$.
If $t\mapsto f(x,t)$ is monotone increasing (resp. decreasing)
then it is necessarily eventually constant.
In that case, the considered notion of approximation from below
(resp. from above) coincides with monotone approximation.
\\
Fix $n\geq1$. An important case is the noetherian set
$(\D,\subseteq)$ of ideals in the ring of $n$-variables polynomials
with real algebraic coefficients
(this last hypothesis insures that $\D$ is countable
with a computable ordering).
%
\subsection{Normalized representations}
\label{ss:normalize}
%
It sometimes proves useful to normalize the $f$ in $\max^\D f$.
\begin{prop}\label{p:normalize}
Let $\bbbx$ be a basic set and $\D=(D,<,\rho)$ be a computable
ordered set.
\\
1. Every $F\in\maxpr[\XD]$ is of the form $F=\max^\D f$ for some
partial computable $f:\bbbx\times\bbbn\to D$,
monotone increasing in its second argument, such that
$\dom(f)=Z\times\bbbn$ where $Z$ is some $\Sigma^0_1$ subset of $D$.
\medskip\\
2. If $F\in\maxpr[\XD]$ has $\Sigma^0_1$ domain then one can suppose
$Z=\dom(F)$.
\end{prop}%
\begin{pf}
1. Let $g:\bbbx\times\bbbn\to D$ be partial computable,
monotone increasing in its second argument, such that
$F=\max^\D g$.
Let $Z=\{{\tt x}:\exists t\ ({\tt x},t)\in \dom(g)\}$ be
the first projection of $\dom(g)$.
Let $\theta:\bbbx\to D$ be the partial computable
function with domain $Z$ such that $\theta({\tt x})$
is the value first obtained in $\{g({\tt s},t):t\in\bbbn\}$
by dovetailing over computations of $g({\tt s},0),g({\tt s},1),...$.
Let also
$$
\Delta_{{\tt x},t}=\{g({\tt x},u): u\leq t\wedge\ g({\tt x},u)
\mbox{ halts in $\leq t$ steps}\}
$$
and define $f$ with domain $Z\times\bbbn$ such that $f({\tt x},t)$
is the greatest element of $\{\theta({\tt x})\}\cup\Delta_{{\tt x},t}$
\medskip\\
2. Observe that $Z$ necessarily contains $\dom(F)$.
If $\dom(F)$ is $\Sigma^0_1$ then
$\widehat{f}=f\!\upharpoonright \!(\dom(F)\times\bbbn)$ is also
partial computable and $\max^\D\widehat{f}=\max^\D f$.
\qed\end{pf}
%
%
\section{Kolmogorov complexities $\kmax[\D]$, $\kmin[\D]$}
\label{s:Kmax}
%
Kolmogorov complexity theory goes through with the $\maxpr[\XD]$
and $\minpr[\XD]$ classes with no difficulty.
\\
First, we recall Kolmogorov complexity over elements of $\D$.
%
\subsection{Kolmogorov complexity $K^D$}
\label{ss:K}
%
Classical Kolmogorov complexity for elements in $D$
is defined as follows
(cf. Kolmogorov, 1965 \cite{kolmo65},
or Li \& Vitanyi \cite{livitanyi},
Downey \& Hirschfeldt \cite{downey},
G\`acs \cite{gacs93}
or Shen \cite{shen}).
\begin{defn}\label{def:Kphi}
Let $\varphi:\words\to D$.
We denote $K_\varphi:D\to\bbbn$ the partial function
with domain $range(\varphi)$ such that
$$
K_\varphi(d) = \min\{|{\tt p}| : \varphi({\tt p})=d\}
$$
I.e., considering words in $\words$ as programs,
$K_\varphi(d)$ is the shortest length of a program
${\tt p}$ mapped onto $d$ by $\varphi$.
\end{defn}
\begin{thm}[Invariance theorem,
Kolmogorov, 1965 \cite{kolmo65}]\label{thm:invariance1}
Let $\bbbx$ be a basic space and
$\D=(D,<,\rho:\bbbn\to D)$ be a computable partially
ordered set.
When $\varphi$ varies in the family $\PR[\words\to D]$
of partial computable functions $\words\to D$,
there is a least $K_\varphi$, up to an additive constant:
$$\exists\varphi\in \PR[\words\to D]\ \
\forall\psi\in \PR[\words\to D]\
\ \ K_\varphi\leqct K_\psi$$
Such $\varphi$'s are said to be optimal in $\PR[\words\to D]$.
\end{thm}
\begin{defn}\label{def:K}
Kolmogorov complexity $K^D:D\to\bbbn$ is
$K_\varphi$ where $\varphi$ is some fixed optimal function
in $\PR[\words\to D]$.
Thus, $K^D$ is defined up to an additive constant.
\end{defn}
Of course, $K^D$ and $K^\bbbn$ are related.
\begin{prop}\label{p:KDKN}
$K^D\circ\rho\eqct K^\bbbn$.
\end{prop}
\begin{pf}
Since
$\PR[\words\to D]=\{\rho\circ\psi:\psi\in\PR[\words\to\bbbn]\}$
and $K_\psi(n)=K_{\rho\circ\psi}(\rho(n))$ for all
$\psi\in\PR[\words\to\bbbn]$,
we see that if $\varphi$ is optimal in $\PR[\words\to\bbbn]$ then
$\rho\circ\varphi$ is optimal in $\PR[\words\to D]$.
\qed\end{pf}
We also observe the following simple fact:
\begin{prop}\label{p:supKD}
$\sup\{K^D(d):d\in X\}=+\infty$ for every
infinite $X\subseteq D$.
\end{prop}
\begin{pf}
The result is well-known for $K^\bbbn$ and it transfers to
$K^D$ using Proposition \ref{p:KDKN}.
\qed\end{pf}
%
%
\subsection{Enumeration theorem for $\maxpr[\XD]$}
\label{ss:enumMaxPr}
%
The classical enumeration theorem for partial computable
functions goes through the $\max$ operator,
leading to an enumeration of $\maxpr[\XD]$.
First, we recall a folklore result on enumeration of {\em monotone}
partial computable functions.
\begin{prop}\label{p:enum}
Let $\bbbx$ be a basic set and $\D=(D,<,\rho)$ be a computable
ordered set.
Let $\PR[\bbbx\times\bbbn\to D,\uparrow]$ be the family of partial
computable functions $\bbbx\times\bbbn\to D$ which are monotone
increasing in their last argument.
There exists a partial computable function
$\psi:\bbbn\times\bbbx\times\bbbn\to D$ such that
$$
\{\psi_n:n\in\bbbn\}=\PR[\bbbx\times\bbbn\to\D,\uparrow]
$$
where $\psi_n:\bbbx\times\bbbn\to D$ denotes the function
$({\tt x},t)\mapsto\psi(n,{\tt x},t)$.
\end{prop}
\begin{pf}
Let $\phi:\bbbn\times\bbbx\times\bbbn\to D$ be a partial computable
function which enumerates the family $\PR[\bbbx\times\bbbn\to D]$
of partial computable functions $\bbbx\times\bbbn\to D$, i.e,
$$
\{\phi_n:n\in\bbbn\}=\PR[\bbbx\times\bbbn\to D]
$$
We modify $\phi$ to $\psi$ so as to get an enumeration of
$\PR[\bbbx\times\bbbn\to D,\uparrow]$.
Consider an injective computable enumeration
$(n_i,{\tt x}_i,t_i,d_i)_{i\in\bbbn}$ of the graph of $\phi$. Let
$$
Z=\{(n_i,{\tt x}_i,t_i,d_i):\forall j<i\
(n_j=n_i\wedge x_j=x_i\wedge t_j<t_i \Rightarrow d_j\leq d_i)\}
$$
Let $\psi:\bbbn\times\bbbx\times\bbbn\to D$ be the partial computable
function with graph $Z$.
It is clear that $\psi$ is monotone increasing in its last argument,
so that so are all $\psi_n$'s. Also, if $\phi_n$ is monotone
increasing in its last argument then
$\{n\}\times{\textit graph}(\phi_n)$
is included in $Z$, so that $\psi_n=\phi_n$.
Thus, the $\psi_n$'s enumerate
$\PR[\bbbx\times\bbbn\to D,\uparrow]$.
\qed\end{pf}
\begin{thm}[Enumeration theorem for {$\maxpr[\XD]$}]
\label{thm:enum}
Let $\bbbx$ be a basic set and $\D=(D,<,\rho)$ be a computable
ordered set.
There exists a function $E:\bbbn\times\bbbx\to D$ in
$\maxpr[\bbbn\times\bbbx\to D]$ such that
$$
\{E_n:n\in\bbbn\}=\maxpr[\XD]
$$
where $E_n:\bbbx\to D$ denotes the function
satisfying $E_n({\tt x})=E(n,{\tt x})$.
\end{thm}
\begin{pf}
Let $\psi:\bbbn\times\bbbx\times\bbbn\to D$ be a partial computable
function which enumerates $\PR[\bbbx\times\bbbn\to D,\uparrow]$.
Let $E:\bbbn\times\bbbx\to D$ be such that
$E_n=\max^\D\psi_n$ for all $n$.
For any $F:\bbbx\to D$ in $\maxpr[\XD]$ there exists $n$ such that
$F=\max^\D\psi_n$. We then have
\begin{eqnarray*}
{\tt x}\in dom(F)
&\Leftrightarrow&
\{\psi_n({\tt x},t):\mbox{$t$ s.t. $\psi_n({\tt x},t)$ is defined}\}
\mbox{ is finite non empty}\\
&\Leftrightarrow&
\{\psi(n,{\tt x},t):\mbox{$t$ s.t. $\psi(n,{\tt x},t)$ is defined}\}
\mbox{ is finite non empty}\\
&\Leftrightarrow&(n,{\tt x})\in dom(E)\\
&\Leftrightarrow&{\tt x}\in dom(E_n)
\\
F({\tt x})&=&\mbox{greatest element of }
\{\psi_n({\tt x},t):\mbox{$t$ s.t. $\psi_n({\tt x},t)$ is defined}\}\\
&=&\mbox{greatest element of }
\{\psi(n,{\tt x},t):\mbox{$t$ s.t. $\psi(n,{\tt x},t)$ is defined}\}\\
&=&E(n,{\tt x})\\
&=&E_n({\tt x})
\end{eqnarray*}
Which proves that $E$ enumerates $\maxpr[\XD]$.
\qed\end{pf}
%
%
\subsection{Kolmogorov complexity $\kmax[\D]$ and $\kmin[\D]$}
\label{ss:Kmax}
%
The invariance theorem extends easily to $\maxpr[\WD]$,
leading to Kolmogorov complexity $\kmax[\D]:D\to\bbbn$.
\begin{thm}[Invariance theorem for {$\maxpr[\WD]$}]
\label{thm:invarianceMaxPr}
Let $\bbbx$ be a basic space and
$\D=(D,<,\rho:\bbbn\to D)$ be a computable partially
ordered set.
When $F$ varies in the family $\maxpr[\WD]$
there is a least $K_F$, up to an additive constant:
$$
\exists U\in \maxpr[\WD]\ \
\forall F\in \maxpr[\WD]\ \ \ K_U\leqct K_F
$$
Such $U$'s are said to be optimal in $\maxpr[\WD]$.
\end{thm}
\begin{pf}
The usual proof works.
Let $E:\bbbn\times\words\to D$ in $\maxpr[\bbbn\times\words\to\D]$
be an enumeration of $\maxpr[\WD]$.
Define $U:\words\to D$ such that $U(0^n1p)=E(n,p)$ and
$U(q)$ is undefined if $q$ is not of the form $0^n1p$ for some
$n\in\bbbn$ and $p\in\words$.
If $F\in\maxpr[\WD]$ and $F=E_n$ then
\begin{eqnarray*}
K_F(d)
&=&\min\{|{\tt p}| : F({\tt p})=d\}\\
&=&\min\{|{\tt p}| : E(n,{\tt p})=d\}\\
&=&\min\{|{\tt p}| : U(0^n1p)=d\}\\
&=&\min\{|0^n1p| : U(0^n1p)=d\}-n-1\\
&\geq&\min\{|{\tt q}| : U({\tt q})=d\}-n-1\\
&=&K_U(d)-n-1
\end{eqnarray*}
\qed\end{pf}
\begin{defn}\label{def:Kmax}
Kolmogorov complexity $\kmax[\D]$ is defined up to an additive
constant as any $K_U$ where $U$ is optimal in $\maxpr[\WD]$.
\\
Kolmogorov complexity $\kmin[\D]$ is $\kmax[\D']$ where
$\D'$ is the reverse order of $\D$.
\end{defn}
%
%
\section{Main theorems:
comparing $K$, $\kmax[\D]$, $\kmin[\D]$, $K^{\emptyset'}$}
\label{s:hierarchy}

%
\subsection{The $\infct$ hierarchy theorem}
\label{ss:hierarchy}
%
The main motivation of this section is to compare the
Kolmogorov complexities
$$
K^D,\ \kmax[\D],\ \kmin[\D],\ K^{\emptyset'}:D\to\bbbn
$$
Comparisons of $\kmax[\D],\kmin[\D]$ and $K^D$ turn out to be
a particular application of more general results
dealing with both $\kmax[\D]$ and $\kmin[\D]$ complexities relative
to two computable orders
$\D_\strong=(D,<_\strong,\rho)$
and $\D_\weak=(D,<_\weak,\rho)$ on the same set $D$,
the strong one $<_\strong$ being an extension
of the weak one $<_\weak$.
A question with naturally arises when considering for instance
the prefix and lexicographic orders on $\Sigma^*$.
\medskip\\
In the case of $\bbbn$ with the natural order or of $\Sigma^*$
with the prefix order, the inequalities
$K^{\emptyset',D}\infct\kmax[\D]\infct K^D$
were obtained (modulo Proposition \ref{p:exitMach})
for the prefix version $H^\infty$
Becher \& Figueira \& Nies \& Picci, 2005 \cite{bechernies},
\\
We state our results as three theorems,
the proofs of which are given in \S\ref{ss:1} to \ref{ss:22}.
\begin{thm}[1st hierarchy theorem]
\label{thm:hierarchy1}
Let $\D=(D,<,\rho)$ be a computable ordered set.
\\
1. $K^{\emptyset',D}\infct\inf(\kmax[\D],\kmin[\D])$
\\
2. $\kmax[\D],\kmin[\D]$ are $\leqct$ smaller than $K^D$
but not simultaneously much smaller:
$$
K^D\leqct(\kmax[\D]+\log(\kmax[\D]))+(\kmin[\D]+\log(\kmin[\D]))
$$
\end{thm}
\begin{thm}[2d hierarchy theorem]\label{thm:hierarchy2}$\\ $
1. If $(D,<)$ contains arbitrarily large finite chains then
$\kmax[\D]$ and $\kmin[\D]$ are $\leqct$ incomparable and both
are $\infct$ smaller than $K^D$.
\\
In fact, a much stronger property holds:
\begin{enumerate}
\item[i.]
$K^D$ is not majorized by a computable function of $\kmin[\D]$,
\item[ii.]
$K^D$ is not majorized by a computable function of $\kmax[\D]$,
\item[iii.]
$\kmax[\D]$ is not majorized by a computable function of $\kmin[\D]$,
\item[iv.]
$\kmin[\D]$ is not majorized by a computable function of $\kmax[\D]$,
\end{enumerate}
I.e., for any total computable function $\alpha:\bbbn\to\bbbn$,
the following sets are infinite
\medskip\\\medskip\centerline{$\begin{array}{ccc}
\{d\in D:\kmax[\D](d)\geq\alpha(\kmin[\D](d))\}
&\ \ ,\ \ &
\{d\in D:K^D(d)\geq\alpha(\kmin[\D](d))\}
\\
\{d\in D:\kmin[\D](d)\geq\alpha(\kmax[\D](d))\}
&\ \ ,\ \ &
\{d\in D:K^D(d)\geq\alpha(\kmax[\D](d))\}
\end{array}$}
2. If $(D,<)$ does not contain arbitrarily large finite chains then
$$
\kmin[\D]\eqct\kmax[\D]\eqct K^D
$$
\end{thm}
\begin{thm}[3d hierarchy theorem]
\label{thm:hierarchy3}
Let $\D_\strong=(D,<_\strong,\rho)$
and $\D_\weak=(D,<_\weak,\rho)$
be two computable orders on the same set $D$
(``$\weak$" and ``$\strong$" stand for ``weak" and ``strong")
such that $<_\strong$ is an extension of $<_\weak$.
\\
1. Let $(*)$ be the following condition
\begin{enumerate}
\item[$(*)$]
For all $k$ there exists a strong chain with $k$ elements
which is a weak antichain.
\end{enumerate}
If $(*)$ holds then
    $\kmax[\D_\strong]\infct\kmax[\D_\weak]$
and $\kmin[\D_\strong]\infct\kmin[\D_\weak]$.
\\
In fact, a much stronger property holds:
$\inf(\kmin[\D_\weak],\kmax[\D_\weak])$ is not majorized by a
computable function of $\kmax[\D_\strong]$ or $\kmin[\D_\strong]$.
I.e., for any total computable function $\alpha:\bbbn\to\bbbn$,
the following sets are infinite
$$
\{d\in D:\inf(\kmin[\D_\weak](d),\kmax[\D_\weak](d))
         \geq\alpha(\kmax[\D_\strong](d))\}
$$
$$
\{d\in D:\inf(\kmin[\D_\weak](d),\kmax[\D_\weak](d))
         \geq\alpha(\kmin[\D_\strong](d))\}
$$
2. Let $(**)$ be the following condition
(which is an effective version, tailored for infinite computations,
of the negation of $(*)$, cf. \S\ref{ss:dilworth}).
\begin{enumerate}
\item[$(**)$]
There exists $k$ such that for every partial computable
$f:\words\times\bbbn\to D$ which is monotone increasing in
its second argument relative to the strong order $<_\strong$
there exist partial computable functions
$f_1,...,f_k:\words\times\bbbn\to D$ which are monotone increasing in
their second argument relative to the weak order $<_\weak$
such that
$$
\{f(p,t):t\in\bbbn\}=\bigcup_{i=1,...,k}\{f_i(p,t):t\in\bbbn\}
$$
\end{enumerate}
If $(**)$ holds then
$\kmin[\D_\strong]\eqct\kmin[\D_\weak]$ and
$\kmax[\D_\strong]\eqct \kmax[\D_\weak]$.
\end{thm}
\begin{cor}\label{cor:hierarchy}
Let $\Sigma$ be a finite or infinite countable alphabet,
let $<_1,<$ be computable orders on $\Sigma$ such that
$<_1$ is partial but non trivial and $<$ is a total extension of
$<_1$.
Consider on $\Sigma^*$ the following orders:
the prefix order $<_\prefix$, the lexicographic orders
$<_{\lexico_1}$ and $<_\lexico$ associated to $<_1$ and $<$.
Then
$$
\kmax[\D_{\prefix}]
\infct  \kmax[\D_{<_{\lexico_1}}]
\infct  \kmax[\D_{<_\lexico}]
\ \ ,\ \
\kmin[\D_{\prefix}]
\infct  \kmin[\D_{<_{\lexico_1}}]
\infct  \kmin[\D_{<_\lexico}]
$$
\end{cor}
\begin{pf}
Let $a,b,c,d\in\Sigma$ be such that
$a<_1b$ and $c<d$ but $c\not<_1d$.
Since $<$ extends $<_1$, $c$ and $d$ are $<_1$ incomparable.
Observe that $\{a^nb:n\in\bbbn\}$ is an infinite increasing chain
for $<_{\lexico_1}$ and an antichain for the prefix order.
Also, $\{c^nd:n\in\bbbn\}$ is an infinite increasing chain for
$<_{\lexico}$ and an antichain for the $<_{\lexico_1}$ order.
This gives condition $(*)$ relative to the pairs
$(<_\prefix,<_{\lexico_1})$ and $(<_{\lexico_1},<_{\lexico})$
of orders on $\Sigma^*$.
\qed\end{pf}

%
\subsection{$(*)$ is an effective version
of the negation of $(**)$}
\label{ss:dilworth}
%
Recall Dilworth's theorem.
\begin{thm}[Dilworth, 1950 \cite{dilworth}]\label{thm:dilworth}
Let $\D=(D,<)$ be an ordered set and $k\in\bbbn$.
If every antichain in $D$ has at most $k$ elements then
$D$ is the union of $k$ chains.
\end{thm}
Dilworth's theorem leads to an equivalent form $(\dagger)$ of $(*)$
and condition $(**)$ appears as an effective version of
$\neg(\dagger)$, tailored for infinite computations.
\begin{prop}\label{p:notstar}
Let $\D_\strong=(D,<_\strong)$ and $\D_\weak=(D,<_\weak)$
be two orders on the same set $D$
such that $<_\strong$ is an extension of $<_\weak$.
\\
Then $(*)$ is equivalent to the following condition $(\dagger)$ :
\begin{enumerate}
\item[$(\dagger)$]
For all $k$ there exists a finite strong chain $X$
which is not the union of $k$ weak chains
\end{enumerate}
\end{prop}
\begin{pf}
$(*)\Rightarrow(\dagger)$.
Apply $(*)$ with $k+1$ and observe that a weak antichain with
$k+1$ elements cannot be the union of $k$ weak chains.
\\
$\neg(*)\Rightarrow\neg(\dagger)$.
Let $k$ be an integer which contradicts $(*)$.
Then, in any strong chain, any weak antichain has $<k$ elements.
Apply Dilworth's theorem to get $\neg(\dagger)$.
\qed\end{pf}
\begin{rem}
1. Clearly $(**)\Rightarrow\neg(\dagger)$.
We do not know whether the converse implication holds or not.
The problem is that the proof of Dilworth's theorem is not
incremental as we now detail.
Let $X\cup\{d\}$ be a strong chain with $d>_\strong x$
for all $x\in X$ and such that every weak antichain included in $X$
has at most $k$ elements.
If $X$ is covered by $k$ weak chains $C_1,...,C_k$ then
$d$ may be incomparable to the top elements of all these $k$ chains.
Thus, though $X\cup\{d\}$ is also the union of $k$ weak chains,
such chains may be quite different from the $C_i$'s.
\\
Condition $(**)$ (as contrasted to $\neg(\dagger)$),
does insure such an incremental character.
\\
2. In case $<_\weak$ has a smallest element $d$, condition $(**)$
is equivalent to the analog condition in which functions
$f,f_1,...,f_k$ are replaced by total computable $g,g_1,...,g_k$.
This can be seen by defining $g,g_1,...,g_k$ from $f,f_1,...,f_k$
as follows
$$
g(p,0)=d\ \ ,\ \
g(p,t+1)=\left\{
\begin{array}{ll}
f(p,t)&\mbox{if $f(p,t)$ converges in $\leq t$ steps}\\
g(p,t)&\mbox{otherwise}
\end{array}\right.
$$
and the same with $g_1,...,g_k$ from $f_1,...,f_k$.
\end{rem}
%
%
\subsection{$\kmax[\D],\kmin[\D]$ are not simultaneously
much smaller than $K^D$}
\label{ss:notsimultaneously}
%
\begin{lem}\label{l:notsimultaneously}
Let $\D=(D,<,\rho)$ be a computable ordered set.
Let $c:\words\times\words\to\words$ be a total computable injective
map and let $J:\bbbn\times\bbbn\to\bbbn$ and $M\in\bbbn$ be such that
$|c(p,q)|\leq J(|p|,|q|)+M$ for all $p,q\in\words$. Then
$$
K^D\leqct J(\kmin[\D],\kmax[\D])
$$
In particular (with the special convention $\log(0)=0$),
$$
K^D\leqct(\kmax[\D]+\log(\kmax[\D]))+(\kmin[\D]+\log(\kmin[\D]))
$$
\end{lem}
\begin{pf}
Let $U,V:\words\to D$ be optimal in $\maxpr[\WD]$ and
$\minpr[\WD]$, i.e. $\kmax[\D]=K_U$ and  $\kmin[\D]=K_V$.
Let $f,g:\words\times\bbbn\to D$ be partial computable, respectively
monotone increasing and decreasing with respect to their 2d argument
such that $U=\max^\D f$ and $V=\min^\D g$.
\medskip\\
Define a partial computable function $\varphi:\words\to D$ as follows:
\begin{itemize}
\item
If $r$ is not in $range(c)$ then $\varphi(r)$ is undefined.
Else, from input $r$, get $p$ and $q$ such that $c(p,q)=r$.
\item
Dovetail computations of the $f(p,t)$'s and $g(q,t)$'s
for $t=0,1,2,...$.
\item
If and when there are $t',t''$ such that $f(p,t')$ and $g(q,t'')$
are both defined and have the same value then
output their common value and halt.
\end{itemize}
By the invariance theorem, there is a constant
$N$ such that $K^D\leq K_\varphi+N$.
\\
Let $d\in D$ and let $p,q$ be shortest programs such that
$U(p)=V(q)=d$, i.e. $\kmax[\D](d)=|p|$ and $\kmin[\D](d)=|q|$.
\\
Observe that, whenever $f(p,t')$ and $g(q,t'')$ are both defined,
we have $f(p,t')\leq d\leq g(p,t'')$.
Also, since
$U(p)=\max^\D\{f(p,t):t\in\bbbn\}$ and
$V(q)=\min^\D\{g(q,t):t\in\bbbn\}$,
there are $t',t''$ such that $f(p,t')=d=g(q,t'')$.
Therefore, $\varphi(c(p,q))$ halts and outputs $d$.
Therefore
$$
K^D(d)\leq K_\varphi(d)+N\leq |c(p,q)|+N
\leq J(\kmin[\D](d),\kmax[\D](d))+N
$$
The last assertion of the Lemma is obtained with the injective map
$$
c(p,q)=\left\{
\begin{array}{ll}
0^{|Bin(|p|)|}1Bin(|p|)pq&\mbox{if }|p|\leq|q|\\
1^{|Bin(|q|)|}0Bin(|q|)pq&\mbox{if }|p|>|q|
\end{array}\right.
$$
(where $Bin(x)$ denotes the binary representation of $x$) since
$$
|c(p,q)|=|p|+|q|+2\lfloor\log(\min(|p|,|q|))\rfloor+3
\leq(|p|+\log(|p|))+(|q|+\log(|q|))+3
$$
\qed\end{pf}
%
%
\subsection{$\kmax[\D],\kmin[\D]$ and the jump}
\label{ss:jump}
%
\begin{prop}\label{p:oracleprime}
1. Let $\bbbx$ be a basic space.
All functions in $\maxpr[\XD]$ and $\minpr[\XD]$
are partial computable in $\emptyset'$.
In particular, $K^D$ is recurcomputablesive in $\emptyset'$.
\medskip\\
2. $\kmin[\D]$ and $\kmax[\D]$ are computable in $\emptyset'$.
\end{prop}
\begin{pf}
1. Proposition \ref{p:syntax} insures that any $F:\XD$ in
$\maxpr[\XD]$ or $\minpr[\XD]$ has $\Sigma^0_1\wedge\Pi^0_1$ graph.
Therefore two calls to oracle $\emptyset'$ suffice to decide
$F({\tt x})=d$.
\medskip\\
2. Let $p_0,p_1,\ldots$ be a {\em length increasing} enumeration
of $\words$ and let $U:\WD$ be optimal in $\maxpr[\WD]$,
i.e. $K_U=\kmax[\D]$.
One can compute $\kmax[\D](d)$ with oracle $\emptyset'$ as follows:
\begin{enumerate}
\item[i.]
Using oracle $\emptyset'$, test successive
equalities $U(p)=d$ (cf. Point 1)
for programs $p = p_0,p_1,\ldots$.
\item[ii.]
When such an equality holds (which necessarily does happen)
then output $|p|$ and halt.
\end{enumerate}
Idem with $\kmin[\D]$.
\qed\end{pf}
%
%
\subsection{Proof of Theorem \ref{thm:hierarchy1}
(1st hierarchy theorem)}
\label{ss:1}
%
1. {\em Large inequality
$K^{\emptyset',D}\leqct\inf(\kmin[\D],\kmax[\D])$.}
Point 1 of Proposition \ref{p:oracleprime} insures that
$\maxpr[\D]$ and $\minpr[\D]$ are included in $\PR[\emptyset']$.
Therefore $K^{\emptyset',D}\leqct\kmin[\D]$ and
       $K^{\emptyset',D}\leqct\kmax[\D]$,
i.e. $K^{\emptyset',D}\leqct\inf(\kmin[\D],\kmax[\D])$.
\\
{\em Strict inequality
$K^{\emptyset',D}\infct\inf(\kmin[\D],\kmax[\D])$.}
Point 2 of Proposition \ref{p:oracleprime} insures that
$\inf(\kmin[\D],\kmax[\D])$ is computable in $\emptyset'$.
Now, the well-known fact that if $\psi\eqct K^D$ then $\psi$ is not
computable relativizes: if $\psi\eqct K^{\emptyset',D}$ then $\psi$
is not computable in $\emptyset'$.
In particular,
$\inf(\kmin[\D],\kmax[\D])\neq_{\rm ct} K^{\emptyset',D}$.
\\
2. This is the contents of Lemma \ref{l:notsimultaneously}.
\hfill{$\Box$}

%
\subsection{Inequalities $\kmax[\D_\strong]\leqct\kmax[\D_\weak]$
and $\kmin[\D_\strong]\leqct\kmin[\D_\weak]$}
\label{ss:easy}
%
The following result is straightforward.
\begin{prop}\label{p:stleqctwk}
With the notations of Theorem \ref{thm:hierarchy3},
$$
\kmin[\D_\strong]\leqct\kmin[\D_\weak]\ \ ,\ \
\kmax[\D_\strong]\leqct \kmax[\D_\weak]
$$
\end{prop}
\begin{pf}
Since $<_\strong$ extends $<_\weak$, every partial computable
function $\words\to D$ which is monotone
increasing in its second argument relative to $<_\weak$ is also
monotone increasing relative to $<_\strong$.
So that $\maxpr[\D_\weak]\subseteq\maxpr[\D_\strong]$.
Which yields $\kmax[\D_\strong]\leqct \kmax[\D_\weak]$.
\qed\end{pf}
%
%
\subsection{If $(*)$ holds: proof of Point 1
            of Theorem \ref{thm:hierarchy3} (3rd hierarchy theorem)}
\label{ss:31}
%
We use the notations of Theorem \ref{thm:hierarchy3}.
\begin{lem}\label{l:star}
Let $\alpha:\bbbn\to\bbbn$ be a total computable function.
\\
If condition $(*)$ holds then there exists total functions
$F,G:\bbbn\to D$ respectively in $\maxr[{\bbbn\to\D}_\strong]$ and
$\minr[{\bbbn\to\D}_\strong]$ and a constant $c$
such that, for all $i\in\bbbn$,
\medskip\\\centerline{$\begin{array}{ccccc}
\kmax[\D_\weak](F(i))\geq\alpha(i)&\ \ ,\ \ &
\kmin[\D_\weak](F(i))\geq\alpha(i)&\ \ ,\ \ &
\kmax[\D_\strong](F(i))\leq \log(i)+c
\\
\kmax[\D_\weak](G(i))\geq\alpha(i)&\ \ ,\ \ &
\kmin[\D_\weak](G(i))\geq\alpha(i)&\ \ ,\ \ &
\kmin[\D_\strong](G(i))\leq \log(i)+c
\end{array}$}
\end{lem}
\begin{pf}
1. Since $(*)$ holds, for all $i\in\bbbn$, there exists a finite
strong chain with $2^{\alpha(i)+1}$ elements which is a weak
antichain.
Dovetailing over subsets of $D$ with $2^{\alpha(i)+1}$ elements,
one can effectively find such a strong chain $Z_i$.
Thus, there exists a total computable function
$\sigma:\bbbn\times\bbbn\to D$ such that, for all $i\in\bbbn$,
\begin{itemize}
\item
$\sigma(i,0)<_\strong\sigma(i,1)<_\strong...
                                <_\strong\sigma(i,2^{\alpha(i)+1}-1)$
\item
$Z_i=\{\sigma(i,j):j=0,...,2^{\alpha(i)+1}-1\}$ is a weak antichain.
\end{itemize}
2. Let $f$ and $g$ are partial computable functions
$\words\times\bbbn\to D$ such that
$U=\max^{\D_\weak}f$ and $V=\min^{\D_\weak}g$ are optimal
in $\maxpr[\WD_\weak]$ and $\minpr[\WD_\weak]$,
i.e. $K_U=\kmax[\D_\weak]$ and $K_V=\kmin[\D_\weak]$.
\\
We observe that inequalities
$\kmax[\D_\weak](F(i))\geq\alpha(i)$ and
$\kmin[\D_\weak](F(i))\geq\alpha(i)$
are equivalent to disequalities $U(p)\neq F(i)$ and $V(p)\neq F(i)$
for every $p$ such that $|p|<\alpha(i)$.
\\
We define $F,G:\bbbn\to D$ as
$F=\max^{\D_\strong}\ell$ and $F=\min^{\D_\strong}\ell$ for some
total computable $\ell:\bbbn\times\bbbn\to D$.
Let
\begin{eqnarray*}
X_p&=&\{f(p,t):\mbox{$t$ s.t. $f(p,t)$ converges}\}\\
Y_p&=&\{g(p,t):\mbox{$t$ s.t. $g(p,t)$ converges}\}
\\
X^t_p&=&
\{f(p,t')\in Z_i
:t'\leq t\mbox{ and $f(p,t')$ converges in $\leq t$ steps}\}
\\
Y^t_p&=&
\{g(p,t')\in Z_i
:t'\leq t\mbox{ and $g(p,t')$ converges in $\leq t$ steps}\}
\end{eqnarray*}
Since $Z_i$ is a weak antichain and $X_p,Y_p$ are weak chains,
each one of the sets $Z_i\cap X_p$ and $Z_i\cap Y_p$
has at most one element.
Thus, $\bigcup_{|p|<\alpha(i)}(X_p\cup Y_p)$ has at most
$2(2^{\alpha(i)}-1)=2^{\alpha(i)+1}-2$
elements in $Z_i$.
Since $Z_i$ has $2^{\alpha(i)+1}$ elements and the $\sigma(i,j)$'s
are in $Z_i$, the following definition makes sense:
$$
\ell(i,t)=
\sigma(i,j)\mbox{ where $j$ is least such that }
\sigma(i,j)\notin\bigcup_{|p|<\alpha(i)}(X^t_p\cup Y^t_p)
$$
Now, $F(i)=(\max^{\D_\strong}\ell)(i)$ and
$G(i)=(\min^{\D_\strong}\ell)(i)$ are of the form $\ell(i,t'_i)$
and $\ell(i,t''_i)$ for some $t'_i,t''_i$,
hence they are not in $\bigcup_{|p|<\alpha(i)}(X_p\cup Y_p)$.
In particular, since
$U(p)={\max}^{\D_\weak} X_p$ is in $X_p$ and
$V(p)={\min}^{\D_\weak} Y_p$ is in $Y_p$,
we see that $F(i)$ and $G(i)$ are not in $\{U(p),V(p)\}$
for any $|p|<\alpha(i)$.
Which proves that
$\kmax[\D_\weak](F(i))$, $\kmin[\D_\weak](F(i))$,
$\kmax[\D_\weak](G(i))$ and $\kmin[\D_\weak](G(i))$ are all
$\geq\alpha(i)$.
\medskip\\
3. Since $F\in\maxr[{\bbbn\to\D}_\strong]$, the invariance theorem
insures that $\kmax[\D_\strong]\leqct K_F$.
Now, $K_F(F(i))\leqct\log(i)$, hence the inequality
$\kmax[\D_\strong](F(i))\leq \log(i)+c$ for some constant $c$.
Idem with $\kmin[\D_\strong](G(i))$.
\qed\end{pf}
\noindent{\bf Proof of Point 1 of Theorem \ref{thm:hierarchy3}}.
Apply Lemma \ref{l:star} with $\alpha'$ such that
$\alpha'$ is monotone increasing and
$\alpha'(i)\geq\max(\alpha(i),i)$ for all $i$.
Since $\alpha'(i)$ tends to $+\infty$ with $i$, so does $F(i)$.
Let $i_0$ be such that $\log(i)+c\leq i$ for all $i\geq i_0$.
Since $\alpha'$ is increasing and $\alpha'\geq\alpha$,
for all $i\geq i_0$ we have
$$
\kmax[\D_\weak](F(i))\geq\alpha'(i)
\geq\alpha'(\lfloor\log(i)+c\rfloor)
\geq\alpha'(\kmax[\D_\strong](F(i)))
\geq\alpha(\kmax[\D_\strong](F(i)))
$$
Similarly, we have
$\kmin[\D_\weak](F(i))\geq\alpha(\kmax[\D_\strong](F(i)))$
and
$\kmax[\D_\weak](G(i))\geq\alpha(\kmin[\D_\strong](G(i)))$
and
$\kmin[\D_\weak](G(i))\geq\alpha(\kmin[\D_\strong](G(i)))$.
\\
Finally, observe that $\{F(i):i\geq i_0\}$ and $\{G(i):i\geq i_0\}$
are infinite.
Which concludes the proof of Point 1 of
Theorem \ref{thm:hierarchy3}.
\hfill{$\Box$}
%
%
\subsection{Proof of Point 1 of Theorem \ref{thm:hierarchy2}
(2d hierarchy theorem)}
\label{ss:21}
%
{\em Comparing $K^D$ to $\kmax[\D]$ and $\kmin[\D]$.}\\
Let $<_\strong$ be $<$ and $<_\weak$ be the empty order.
Then
$$
\kmax[\D_\strong]=\kmax[\D]\ \ , \ \
\kmin[\D_\strong]=\kmin[\D]\ \ , \ \
\kmax[\D_\weak]=\kmin[\D_\weak]=K^D
$$
The condition (in Point 1 of Theorem \ref{thm:hierarchy2})
that $\D$ contains arbitrarily large chains insures condition
$(*)$ about $<_\strong$ and $<_\weak$.
Thus, we can apply (the just proved) Point 1 of Theorem
\ref{thm:hierarchy3}. This gives properties i and ii of Point 1
of Theorem \ref{thm:hierarchy2}.
\medskip\\
{\em Comparing $\kmax[\D]$ and $\kmin[\D]$.}\\
We shall prove properties iii and iv of Point 1
of Theorem \ref{thm:hierarchy2} using properties i and ii
and also Lemma \ref{l:notsimultaneously}.
\\
Applying Lemma \ref{l:notsimultaneously}, let $c$ be such that,
$$
(\dagger)\hspace{1cm}K^D\leq 2\,(\kmax[\D]+\kmin[\D])+c
$$
Property iii applied to $\alpha'(i)=2\,(\alpha(i)+i)+c$ insures that
the set
$$
X=\{d:K^D(d)\geq2\,(\alpha(\kmin[\D](d))+\kmin[\D](d))+c\}
$$
is infinite.
Now, using $(\dagger)$, we see that, for $d\in X$,
$$
2\,(\alpha(\kmin[\D](d))+\kmin[\D](d))+c\leq K^D(d)
\leq2\,(\kmax[\D](d)+\kmin[\D](d))+c
$$
hence $\kmax[\D](d)\geq\alpha(\kmin[\D](d))$.
Which proves iii.
The proof of iv is similar.
\hfill{$\Box$}
%
%
\subsection{If $(**)$ holds: proof of Point 2
of Theorem \ref{thm:hierarchy3} (3d hierarchy theorem)}
\label{ss:32}
%
\begin{lem}\label{l:starstar}
With the notations of Theorem \ref{thm:hierarchy3},
if condition $(**)$ holds then
$$
\kmin[\D_\strong]\geqct\kmin[\D_\weak]\ \ ,\ \
\kmax[\D_\strong]\geqct \kmax[\D_\weak]
$$
\end{lem}
\begin{pf}
1. Let $k$ be as in $(**)$.
Let $U_\strong$ be optimal in $\maxpr[\D_\strong]$ and
$f:\words\times\bbbn\to D$ be partial computable such that
$\max^{\D_\strong}f=U_\strong$.
\\
Due to Proposition \ref{p:normalize}, we can suppose that $f$ has
domain of the form $Z\times\bbbn$ and is monotone increasing in its
second argument, with respect to the strong order.
\\
Applying $(**)$ to $f$, we get $k$ partial computable functions
$f_1,...,f_k$, monotone increasing in their second argument,
with respect to the weak order, such that
$$
(\sharp)\hspace{1cm}
\{f(p,t):t\in\bbbn\}=\bigcup_{i=1,...,k}\{f_i(p,t):t\in\bbbn\}
$$
Define $g:\words\times\bbbn\to D$ such that
$$
g(q,t)=\left\{
\begin{array}{ll}
f_i(p,t)&\mbox{if $q=0^i1^{k-i}p$ for some $p$ and $1\leq i\leq k$}\\
\mbox{undefined}&\mbox{otherwise}
\end{array}\right.
$$
Clearly, $g$ is partial computable and monotone increasing
in its second argument relative to the weak order $<_\weak$.
\\
If $p\in dom(U_\strong)$, then $\{f(p,t):t\in\bbbn\}$ is finite and
non empty. Let $f(p,t_p)$ be its $<_\strong$ greatest element.
Condition $(\sharp)$ insures that there exists $i$ such that
$\{g(0^i1^{k-i}p,t):t\in\bbbn\}$ is finite and contains $f(p,t_p)$.
Since $g$ is $<_\weak$ increasing in $t$, the set
$\{g(0^i1^{k-i}p,t):t\in\bbbn\}$ is a weak chain.
Since $<_\strong$ extends $<_\weak$, $f(p,t_p)$ is necessarily
its $<_\weak$ greatest element.
Thus,
$$
U_\strong(p)=f(p,t_p)=({\max}^{\weak} g)(0^i1^{k-i}p)
$$
This proves that, for all $d\in D$,
\begin{eqnarray*}
\kmax[D_\strong](d)&=&\mbox{least $|p|$ such that $U_\strong(p)=d$}
\\ &=&\mbox{least $|p|$ such that
            $({\max}^{\weak} g)(0^i1^{k-i}p)=d$ for some $i$}\\
&\geq&\mbox{least $|q|-k$ such that $({\max}^{\weak} g)(q)=d$}\\
&=&K_{\max^{\D_\weak}g}(d)-k
\end{eqnarray*}
Since, by the invariance theorem,
$K_{\max^{\D_\weak}g}\geqct\kmax[D_\weak]$,
we get the desired inequality
$\kmax[D_\strong]\geqct\kmax[D_\weak]$.
\\
2. Considering the reverse orders, we get the inequality
$\kmin[\D_\strong]\geqct\kmin[\D_\weak]$.
\qed\end{pf}
\noindent{\bf Proof of Point 2 of Theorem \ref{thm:hierarchy3}}.
Straightforward from the above Lemma \ref{l:starstar} and
Proposition \ref{p:stleqctwk}.
\hfill$\Box$
%
%
\subsection{Proof of Point 2 of Theorem \ref{thm:hierarchy2}
(2d hierarchy theorem)}
\label{ss:22}
%
As in \S\ref{ss:21}, let
$<_\strong$ be $<$ and $<_\weak$ be $\emptyset$, so that
$$
\kmax[\D_\strong]=\kmax[\D]\ \ , \ \
\kmin[\D_\strong]=\kmin[\D]\ \ , \ \
\kmax[\D_\weak]=\kmin[\D_\weak]=K^D
$$
Suppose all chains in $(D,<)$ have length $\leq k$.
We shall prove condition $(**)$ for the above orders
$<_\strong$ and $<_\weak$.
\\
Let $f:\words\times\bbbn\to D$ be partial computable,
monotone increasing in its 2d argument for the strong order,
i.e. for the $<$ order.
Compute $f(p,t)$ for $t=0,1,...$ to get the $\leq k$ distinct
elements of the chain $\{f(p,t):t\in\bbbn\}$ (not necessarily in
increasing order) and let $f_i(p)$ be the $i$-th element so
obtained (if there is some).
Then $f_0,...,f_k:\words\to D$ are partial computable and
$$
\{f(p,t):t\in\bbbn\}
=\{f_i(p):\mbox{ $i$ s.t. $f_i(p)$ is defined}\}
$$
which insures condition $(**)$.
\\ Applying Point 2 of Theorem \ref{thm:hierarchy3} (proved above),
we get $\eqct$ equalities which are exactly those of Point 2 of
Theorem \ref{thm:hierarchy2}.
\hfill$\Box$
%
%
%
%
%
\section{Complementary results about the $Max$ and $Min$ classes}
\label{s:complements}
%
In this section we further investigate the different
$Max$ and $Min$ classes.
The results do not involve as many technicalities as those
of \S\ref{s:hierarchy}.
%
%
\subsection{Total functions in $\maxr[\XD]$ and $\maxpr[\XD]$}
\label{ss:total}
%
As a straightforward corollary of Point 2 of
Proposition \ref{p:normalize}, we get the following result.
\begin{thm}\label{thm:MaxPrRec}
The classes $\maxr[\XD]$ and $\maxpr[\XD]$
contain the same total functions:
$$
\maxpr[\XD]\cap D^\bbbx=\maxr[\XD]\cap D^\bbbx
$$
\end{thm}
%
%
\subsection{Comparing $\maxpr[\XD],\maxr[\XD]$ and
$\PR[\XD],\Rec[\XD]$}
\label{ss:PRversusMaxPR}
%
\begin{prop}\label{p:PRversusMaxPR}
Let $\bbbx$ be a basic set and $\D=(D,<,\rho)$ be a computable
ordered set.
\\
1. If $<$ is empty then $\PR[\XD]=\maxpr[\XD]$ and
$\Rec[\XD]=\maxr[\XD]$.
\medskip\\
2. If $<$ is not empty then $\maxr[\XD]$ contains non computable
total functions.
In particular, $\PR[\XD]\subset\maxpr[\XD]$ and
$\Rec[\XD]\subset\maxr[\XD]$
(where $\subset$ denotes strict inclusion).
\medskip\\
3. Whatever be $<$,
$\PR[\XD]$ is not included in $\minr[\XD]\cup\maxr[\XD]$.
\end{prop}
\begin{pf}
1. Straightforward.
\\
2. Inclusions $\PR[\XD]\subseteq\maxpr[\XD]$ and
$\Rec[\XD]\subseteq\maxr[\XD]$ are obvious.
\\
Suppose there exists comparable distinct elements $a<b$ in $D$.
Let $Z$ be some computably enumerable non computable subset
of $\bbbx$ and let $\theta:\bbbn\to\bbbx$ be a total computable
map with range $Z$.
Define $f:\bbbx\times\bbbn\to D$ total computable,
monotone increasing in $t$, such that
$$
f({\tt x},t)=\left\{
\begin{array}{ll}
a&\mbox{if }{\tt x}\notin\{\theta(n):n\leq t\}
\\
b&\mbox{otherwise}
\end{array}\right.
$$
Then $\max f$ is total and $(\max f)^{-1}(b)=Z$
and $(\max f)^{-1}(a)=\bbbx\setminus Z$.
Since $Z$ is not computable, $\max f$ is not
computable. Which proves $\Rec[\XD]\subset\maxr[\XD]$
\\
3. First, we consider the case where $(D,<)$ has a minimal element $d$.
Let $\pi^Z_d:\bbbx\to\ D$ be the partial computable function
with domain $Z$ (as in Point 2 of this proof) which is constant
on $Z$ with value $d$. We show that $\pi^Z_d$ is not in $\maxr[\D]$.
Suppose $f:\bbbx\times\bbbn\to D$ is total computable,
monotone in its second argument, such that $\max^\D f=\pi^Z_d$.
Since $d$ is minimal in $D$,
$(\max^\D f)({\tt x})=d$ if and only if $\forall t\ f({\tt x},t)=d$.
Thus, the computably enumerable set $Z$ would be $\Pi^0_1$,
hence computable,
contradiction.
\\
We now consider the case where $(D,<)$ has no minimal element.
Let $\gamma:D\to D$ be the total computable function which associates
to each $d\in D$ the element $\rho(k_d)$ where $n_d$ is the least
$k$ such that $\rho(k)<d$.
Let $(\phi)_{{\tt e}\in\bbbx}$ be an enumeration of
$\PR[\bbbx\times\bbbn\to D]$ which is
partial computable as a function
$\Phi:\bbbx\times\bbbx\times\bbbn\to D$.
We consider an enumeration $({\tt e}_n,{\tt x}_n,t_n,d_n)_{n\in\bbbn}$ of
the graph of $\Phi$ and define a partial computable function
$\varphi:\XD$ as follows:
$$
\varphi({\tt x})=\left\{
\begin{array}{ll}
\gamma(d_n)&\mbox{if $n$ is least such that ${\tt e}_n={\tt x}_n={\tt x}$}\\
\mbox{undefined}&\mbox{if there is no such $n$}
\end{array}\right.
$$
It is clear that, for every ${\tt e}$,
if $\phi_{\tt e}({\tt e},t)$ is defined for some $t$
then $\varphi({\tt e})$ is defined and
$\varphi({\tt e})<\phi_{\tt e}({\tt e},t)$.
In particular, if $\phi_{\tt e}$ is total then
$\varphi({\tt e})<(\max^\D\phi_{\tt e})({\tt e})$, hence
$\varphi\neq\max^\D\phi_{\tt e}$.
Which proves that $\varphi$ is not in $\maxr[\XD]$.
\\
Arguing with $\Dd$ we get some function in $\PR[\XD]$ which is not
in $\minr[\XD]$.
Considering $\varphi_0,\varphi_1\in\PR[\XD]$ such that
$\varphi_0\notin\maxr[\XD]$ and $\varphi_1\notin\minr[\XD]$
and a computable bijection $\sigma:\bbbx\times\{0,1\}\to\bbbx$
we get a partial computable function $\varphi:\XD$ which is not
in $\maxr[\XD]\cup\minr[\XD]$ by setting
$\varphi(\sigma({\tt x},0))=\varphi_0({\tt x})$ and
$\varphi(\sigma({\tt x},1))=\varphi_1({\tt x})$.
\qed\end{pf}
%
%
\subsection{Post hierarchy and the $Max/Min$ classes}
\label{ss:syntax}
%
We keep notations of \S\ref{ss:domains}.
\begin{thm}\label{thm:syntax}
Let $\bbbx$ be a basic set and $\D$ be a computable ordered set.
\\
1. Let $D'$ be an initial segment of $D$
(i.e. $d'\in D'\wedge e<d'\Rightarrow e\in D'$).
Suppose $D'$ is $\Pi^0_1$ and does not contain any strictly increasing
infinite sequence $d'_0<d'_1<...$. Then
\\
i. Every $D'$-valued function in $\maxpr[\XD]$ has
$\Sigma^0_1\wedge\Pi^0_1$ domain.
\\
ii. Every $D'$-valued function in $\maxr[\XD]$ has $\Pi^0_1$ domain.
\medskip\\
2. Let $D'$ be a final segment of $D$
(i.e. $d'\in D'\wedge e>d'\Rightarrow e\in D'$).
Suppose $D'$ is $\Sigma^0_1$ and does not contain any strictly
increasing infinite sequence. Then
i. Every $D'$-valued function in $\maxpr[\XD]$ has $\Sigma^0_1$ domain.
\\
ii. Every $D'$-valued function in $\maxr[\XD]$ is total.
\end{thm}
\begin{pf}
1. Suppose that $\max^\D f$ is $D'$-valued.
Since $D'$ is an initial segment
and $f$ can be supposed monotone increasing in its second argument,
if $(\max^\D f)({\tt x})$ is defined then, for all $t$,
$f({\tt x},t)$ is either undefined or in $D'$.
Now, since $D'$ has no infinite increasing sequence, the set
$\{f({\tt x},t):t\in\bbbn\mbox{ s.t. }f({\tt x},t)\in D'\}$
cannot be infinite. Thus,
${\tt x}\in \dom(\max^\D f)$ if and only if
$$
\exists t\ f({\tt x},t)\mbox{ is defined } \wedge\
\forall t\ (f({\tt x},t)\mbox{ is defined }\Rightarrow\
f({\tt x},t)\in D')
$$
In case $f$ is total computable, then the above equivalence
is simply
$${\tt x}\in \dom({\max}^\D f)\ \Leftrightarrow\
\forall t\ f({\tt x},t)\in D'$$
2. Since $D'$ is a final segment
and $f$ can be supposed monotone increasing in its second argument,
if $(\max^\D f)({\tt x})$ is defined then, for all $t$ large enough,
$f({\tt x},t)$ is either undefined or in $D'$.
Now, since $D'$ has no infinite increasing sequence, the set
$\{f({\tt x},t):t\in\bbbn\mbox{ s.t. }f({\tt x},t)\in D'\}$
cannot be infinite. Thus,
\begin{eqnarray*}
{\tt x}\in \dom({\max}^\D f)&\Leftrightarrow&
\exists t\ (f({\tt x},t)\mbox{ is defined } \wedge f({\tt x},t)\in D')
\end{eqnarray*}
\qed\end{pf}
The next corollary is an application of the above theorem with
the reverse of the following $\D$'s:
\begin{itemize}
\item
$\D$ is the natural order on $\bbbz$ and $D'=\bbbn$,
\item
$\D$ is the natural order on $\bbbn$
or of the prefix order on $\Sigma^*$ and $D'=D$,
\end{itemize}
\begin{cor}\label{cor:Syntax}
1. Every $\bbbn$-valued function in
$\minpr[\bbbx\to\bbbz]$ (resp. $\minr[\bbbx\to\bbbz]$)
has $\Sigma^0_1\wedge\Pi^0_1$ (resp. $\Pi^0_1$) domain.
\\
2. Let $\D$ be $\bbbn$ with the natural order or $\Sigma^*$
with the prefix partial order.
Then every function
in $\minpr[\bbbx\to\D]$ (resp. $\minr[\bbbx\to\D]$)
has $\Sigma^0_1$ domain (resp. is total).
\end{cor}
%
%
\subsection{$Max\cap Min$ classes}
\label{ss:MaxInterMin}
%
\begin{thm}\label{thm:MaxInterMin}
Let $\bbbx$ be a basic set and $\D=(D,<,\rho)$ be a computable
ordered set.
\\
1. Every function $F:\bbbx\to D$ in $\maxpr[\XD]\cap\minpr[\XD]$
is the restriction of a partial computable function $\XD$ to some
$\Sigma^0_1\wedge\Pi^0_1$ subset of $\bbbx$.
\\
In particular, every total function in
$\maxpr[\XD]\cap\minpr[\XD]$ is computable.
\medskip\\
2. Suppose $\D$ has no maximal (resp. minimal) element.
Then the restriction of any partial computable function $\XD$ to any
$\Sigma^0_1\wedge\Pi^0_1$ subset of $\bbbx$ is in $\maxpr[\XD]$
(resp. $\minpr[\XD]$).
\medskip\\
3. Suppose $\D$ has no maximal or minimal element.
Then $\maxpr[\XD]\cap\minpr[\XD]$ coincides with the family of
restrictions of partial computable functions $\XD$ to
$\Sigma^0_1\wedge\Pi^0_1$ subsets of $\bbbx$.
\end{thm}
\begin{pf}
1. Let $F=\max^\D f=\min^\D g$ where $f,g:\bbbx\times\bbbn\to D$ are
partial computable and $f$ (resp. $g$) is monotone increasing
(resp. decreasing) in its second argument.
Let's check that $F({\tt x})$ is defined if and only if
$$
(*)\hspace{1cm}
(\exists t',t''\ f({\tt x},t')=g({\tt x},t''))
\wedge
(\forall u,v\ f({\tt x},u)\leq g({\tt x},v))
$$
In fact, if $F({\tt x})$ is defined then
\begin{enumerate}
\item[]
$F({\tt x})=f({\tt x},t')=g({\tt x},t'')$ for some $t',t''$,
\item[]
$g({\tt x},u)\leq F({\tt x})\leq f({\tt x},v)$
for all $u,v$ such that $g({\tt x},u),f({\tt x},v)$ are defined.
\end{enumerate}
Conversely, from $(*)$ we see that, for $u\geq t'$ and $v\geq t''$,
$f({\tt x},u)=f({\tt x},t')=g({\tt x},t'')=g({\tt x},v)$.
Hence the finiteness of $\{f({\tt x},u):u\}$ and $\{g({\tt x},v):v\}$.
\\
This proves that the domain of $F$ is $\Sigma^0_1\wedge\Pi^0_1$.
\\
Let $G:\XD$ be the partial computable function defined as follows:
\begin{quote}
Dovetail computations of
$f({\tt x},0),f({\tt x},1),\ldots,
g({\tt x},0),g({\tt x},1),\ldots$
until we get $t',t''$ such that $f({\tt x},t'),g({\tt x},t'')$
are both defined and equal.
Output this common value.
\end{quote}
Applying $(*)$, if $F({\tt x})$ is defined, then so is $G({\tt x})$
and $F({\tt x})=G({\tt x})$.
Thus, $F$ is the restriction of a partial computable function to
some $\Sigma^0_1\wedge\Pi^0_1$ set.
\medskip\\
2. Suppose there is no maximal element.
Since the order $<$ is computable, by dovetailing, one
can define a total computable function $\gamma:D\to D$ such that
$\gamma(d)>d$ for all $d\in D$.
Let $F:\XD$ be partial computable and let $Z\subseteq\bbbx$ be
$\Sigma^0_1\wedge\Pi^0_1$ definable:
$$
{\tt x}\in Z\Leftrightarrow
(\exists t\  R({\tt x},t))\wedge(\forall t\ S({\tt x},t))
$$
where $R,S\subseteq\bbbx\times\bbbn$ are computable.
Letting $\gamma^{(t)}$ denote the $t$-th iterate of $\gamma$,
we define $f:\bbbx\times\bbbn\to D$ as follows:
\begin{eqnarray*}
f({\tt x},t)&=&\left\{
\begin{array}{ll}
F({\tt x})&\mbox{if $F({\tt x})$ converges in $\leq t$ steps}
\\ &\mbox{and }
(\exists t'\leq t\ R({\tt x},t'))\wedge(\forall t'\leq t\ S({\tt x},t'))
\\
\gamma^{(t)}(F({\tt x}))
&\mbox{if $F({\tt x})$ converges in $\leq t$ steps}
\\ &\mbox{and }\exists t'\leq t\ \neg S({\tt x},t')
\\
\mbox{undefined}&\mbox{otherwise}
\end{array}
\right.
\end{eqnarray*}
It is easy to check that $\max^\D f$ is the restriction of $F$
to $Z$.
\\
The assertion with $\minpr[\XD]$ is obtained with the order reverse
to $\D$.
\\
3. Straightforward from Points 1 and 2.
\qed\end{pf}
\begin{rem}
Theorem \ref{thm:syntax} shows that Points 2, 3 of the above
theorem do not hold for general ordered sets $\D$.
\end{rem}
%
%
\section{$\maxr[\WD]$ and $\minr[\WD]$ and Kolmogorov complexity}
\label{s:KwithMaxRec}
%
Since there is no computable enumeration of total computable functions,
it seems a priori desperate to get an invariance theorem for
the class $\maxr[\XD]$.
Nevertheless, there are important cases where such a result
does hold.
For instance, when $\D$ is $\bbbn$ with its usual ordering.
\\
The purpose of this section is to characterize the orders $\D$
such that an invariance theorem holds for the class
$\maxr[\WD]$ (resp. $\minr[\WD]$).
\\
First, we deal with the enumeration theorem.
%
\subsection{ $\maxr[\XD]$ and the enumeration theorem}
\label{ss:enumMaxRec}
%
\begin{thm}[Enumeration theorem for {$\maxr[\XD]$}]
\label{thm:enumMaxr}
Let $\bbbx$ be a basic set and $\D=(D,<,\rho)$ be a computable
ordered set.
The following conditions are equivalent:
\\
i. There exists a smallest element in $\D$.
\\
ii. There exists a function $\widetilde{E}:\bbbn\times\bbbx\to D$
in $\maxr[\bbbn\times\bbbx\to D]$ such that
$$
\{\widetilde{E}_n:n\in\bbbn\}=\maxr[\XD]
$$
where $\widetilde{E}_n:\bbbx\to D$ denotes the function
${\tt x}\mapsto\widetilde{E}(n,{\tt x})$.
\end{thm}
\begin{pf}
$i\Rightarrow ii$.
Let $\alpha\in D$ be the smallest element of $D$.
As in \S\ref{ss:enumMaxPr}, let
$\psi:\bbbn\times\bbbx\times\bbbn\to D$ be partial computable
monotone increasing in its last argument such that
$E=\max^\D\psi$ is an enumeration of $\maxpr[\XD]$.
Consider an injective computable enumeration
$(n_i,{\tt x}_i,t_i,d_i)_{i\in\bbbn}$ of the graph of $\psi$.
Since $\alpha$ is the smallest element, we can define a {\em total}
computable function
$\widetilde{\psi}:\bbbn\times\bbbx\times\bbbn\to D$
as follows:
\begin{eqnarray*}
X(n,{\tt x},t)&=&
\{d_i:i\leq t\wedge n_i=n\wedge {\tt x}_i={\tt x}\wedge t_i\leq t\}
\\
\widetilde{\psi}(n,{\tt x},t)&=&\mbox{greatest element of }
\{\alpha\}\cup X(n,{\tt x},t)
\end{eqnarray*}
Suppose $\psi_n$ is total, we show that
$\max^\D\widetilde{\psi}_n=\max^\D\psi_n$.
Fix some ${\tt x}$.
Observe that
$\{\widetilde{\psi}_n({\tt x},t):t\in\bbbn\}$
is $\{\psi_n({\tt x},t):t\in\bbbn\}$ or
$\{\alpha\}\cup\{\psi_n({\tt x},t):t\in\bbbn\}$.
Thus, $\{\widetilde{\psi}_n({\tt x},t):t\in\bbbn\}$ and
$\{\psi_n({\tt x},t):t\in\bbbn\}$ are simultaneously
finite or infinite, and when finite they have the same greatest
element. Since $\psi_n$ is total, this proves that
$(\max^\D\widetilde{\psi}_n)({\tt x})=(\max^\D\psi_n)({\tt x})$.
Thus, every function in $\maxr[\bbbx\to D]$ is of the form
$\max^\D\widetilde{\psi}_n$ for some $n$.
\\
Set $\widetilde{E}=\max^\D\widetilde{\psi}$. Then
$\widetilde{E}$ is in $\maxr[\bbbn\times\bbbx\to D]$ and
the $\widetilde{E}_n$'s enumerate $\maxr[\bbbx\to D]$.
\medskip\\
$ii\Rightarrow i$. We prove $\neg i\Rightarrow\neg ii$.
Suppose $\D$ has no minimum element.
By dovetailing one can define a total computable map
$\gamma:D\to D$ such that
$d\not\leq\gamma(d)$ for all $d$.
\\
Let $E=\max^\D g:\bbbn\times\bbbx\to D$ where
$g:\bbbn\times\bbbx\times\bbbn\to D$ is total computable monotone
increasing in its last argument.
We define a total computable map $f:\bbbx\to D$ such that
$f\neq E_n$ for all $n$.
Let $\theta:\bbbn\to\bbbx$ be some computable bijection.
Set $f(\theta(n))=\gamma(g(n,\theta(n),0))$.
Then
$$
g(n,\theta(n),0)\not\leq f(\theta(n)) \mbox{ and }
g(n,\theta(n),0)\leq({\max}^\D g)(\theta(n))=E_n(\theta(n))
$$
Thus, $f(\theta(n))\neq E_n(\theta(n))$.
Hence $f\neq E_n$ for all $n$.
\qed\end{pf}
%
%
\subsection{$\maxr[\WD]$ and the invariance theorem}
\label{ss:KmaxRec}
%
If $\D$ contains a smallest element then the enumeration theorem
of \S\ref{ss:enumMaxRec} allows to get an invariance result for
the class $\maxr[\WD]$.
\\
Surprisingly, it turns out that an invariance result can be proved
for partially ordered sets with no smallest element, hence which
fail the enumeration theorem.
\\
Also, in case the class $\maxr[\WD]$ has optimal functions then
they prove to be also optimal for the bigger class $\maxpr[\WD]$.
\begin{thm}\label{thm:KmaxRec}
Let $\bbbx$ be a basic space and
$\D=(D,<,\rho:\bbbn\to D)$ be a computable partially
ordered set.
Let $(*)$ be the following condition on $\D$ :
\begin{enumerate}
\item[$(*)$]
The set of minimal elements of $D$ is finite
and every element of $D$ dominates a minimal element
\end{enumerate}
1. If $\D$ satisfies $(*)$ then
\begin{enumerate}
\item[i.]
Every function in $\maxpr[\WD]$ has an extension
(not necessarily total) in $\maxr[\WD]$.
\item[ii.]
The invariance theorem holds for $\maxr[\WD]$.
\item[iii.]
Every $U$ in $\maxr[\WD]$ which is optimal for $\maxr[\WD]$
is also optimal for the class $\maxpr[\WD]$.
\\
In particular, the Kolmogorov complexity associated to $\maxr[\WD]$
coincides (up to a constant) with that associated to $\maxpr[\WD]$.
\end{enumerate}
2. If $\D$ does not satisfy $(*)$ then the invariance theorem fails
for $\maxr[\WD]$. Moreover, counterexamples can be taken in
the class $\Rec[\WD]$ of total computable functions $\WD$ :
$$
\forall G\in \maxr[\WD]\ \
\exists F\in \Rec[\WD]\ \ \ K_G\not\leqct K_F
$$
\end{thm}
\begin{pf}
1. Suppose $(*)$ holds and let $M=\{m_0,...,m_k\}$ be the
set of minimal elements.
For $i\leq k$, let $D_i=\{d\in D:d\geq m_i\}$.
Some of the $D_i$'s may be finite, though not all of them
(else $D$ would be finite).
Let $\ell\leq k$ be such that $D_i$ is infinite for $i\leq\ell$
and finite for $\ell<i\leq k$.
Since the $D_i$ are computable, for $i\leq\ell$, there exists a
computable map $\rho_i:\bbbn\to D_i$ such that
$\D_i=(D_i,<\cap\,(D_i\times D_i),\rho_i)$ is a computable partially
ordered set.
\\
A. Since $D_i$ has a smallest element, namely $m_i$,
$\maxr[\D_i]$ satisfies the enumeration theorem
(cf. Theorem \ref{thm:enumMaxr}). The proof of Theorem
\ref{thm:invarianceMaxPr} applies, insuring that $\maxr[\D_i]$
satisfies the invariance theorem.
\\
Let $g_i:\words\times\bbbn\to D_i$ be total computable such that
$\max^{\D_i}g_i=U_i:\words\to D_i$ is optimal in $\maxr[\D_i]$.
\\
Let's check that $U_i$ is also optimal in $\maxpr[\D_i]$.
Let $F_i\in\maxpr[\D_i]$ and $F_i=\max^{\D_i}f_i$ where
$f_i:\words\times\bbbn\to D_i$ is partial computable
monotone increasing in its second argument and has domain
$Z_i\times\bbbn$ where $Z_i\subseteq\words$ is computably
enumerable (cf. Proposition \ref{p:normalize}).
Define a total computable map
$\widetilde{f_i}:\words\times\bbbn\to D_i$
such that
$$\widetilde{f_i}(p,t)=\left\{
\begin{array}{ll}
f_i(p,t)&\mbox{if $p$ is seen to be in $Z_i$ in $\leq t$ steps}\\
m_i&\mbox{otherwise}
\end{array}\right.
$$
Set $\widetilde{F_i}=\max^{\D_i}\widetilde{f_i}$.
If $p\in Z$ then $\widetilde{f_i}(p,t)=f_i(p,t)$ for $t$ large
enough, so that
$F_i(p)=(\max^{\D_i}f_i)(p)
=(\max^{\D_i}\widetilde{f_i})(p)=\widetilde{F_i}(p)$.
Thus, $\widetilde{F_i}$ extends $F_i$. Which trivially yields
$K_{\widetilde{F_i}}\leq K_{F_i}$.
Since $\widetilde{F_i}\in\maxr[\D_i]$, we have
$K_{U_i}\leqct K_{\widetilde{F_i}}$. Hence $K_{U_i}\leqct K_{F_i}$.
\\
B. We group the functions $g_i$ and $U_i$ of Point A to get a
total computable $g:\words\times\bbbn\to D$ and the associated
$U=\max^\D g$ in $\maxr[\D]$.
Define $g$ as follows:
$$
g(q,t)=\left\{
\begin{array}{ll}
g_i(p,t)&\mbox{ if $q$ is of the form $0^i1p$ with $i\leq\ell$,
$p\in\words$}
\\
m_0&\mbox{otherwise}
\end{array}\right.
$$
For $i\leq\ell$ and $d\in D_i$, we have
$$
K_U(d)\leq K_{U_i}(d)+i+1
\mbox{ for all $i\leq\ell$ and $d\in D_i$}
$$
Suppose $F$ is in $\maxpr[\D]$ is of the form $F=\max^\D f$ where
$f:\words\times\bbbn\to D$ is partial computable.
For $i\leq\ell$, let $F_i=\max^{\D_i} f_i$ where $f_i:\words\to D_i$
is such that
$$
f_i(p,t)=\left\{
\begin{array}{ll}
f(p,t)&\mbox{ if $f(p,t)$ is defined and is in $D_i$}
\\
\mbox{undefined}&\mbox{otherwise}
\end{array}\right.
$$
Clearly, $F_i$ is the restriction of $F$ to $F^{-1}(D_i)$.
Thus, $K_F(d)=K_{F_i}(d)$ for all $d\in D_i$.
\\
Since $F_i\in\maxpr[\D_i]$ and $U_i$ is optimal in $\maxpr[\D_i]$,
there exists $c_i$ such that $K_{U_i}\leq K_{F_i}+c_i$.
Thus, for $d\in D_i$, we have
$$
K_U(d)\leq K_{U_i}(d)+i+1\leq K_{F_i}(d)+c_i+i+1\leq K_F(d)+c_i+i+1
$$
Let $a$ be the maximum value of $K_F$ on the finite set
$\bigcup_{\ell<j\leq k}D_j$.
Set $c=\sup(\{c_i+i+1:i\leq\ell\}\cup\{a\})$.
Then $K_U(d)\leq K_F(d)+c$ for all $d\in D$.
Which proves that $U$, which is in $\maxr[\D]$, is optimal in
$\maxpr[\D]$.
\\
C. If $V$ in $\maxr[\WD]$ is optimal for
$\maxr[\WD]$ then $K_V\leqct K_U$ (where $U$ is as in B).
Since $U$ is is optimal in $\maxpr[\D]$, so is $V$.
\medskip\\
2. Suppose $(*)$ fails.
Observe that, for every finite subset $Z$ of $D$,
there exists $d$ such that $z\not\leq d$ for all $z\in Z$.
Else, the set of minimal elements of $Z$ would satisfy $(*)$.
\\
Let $D^{<\omega}$ be the set of finite sequences of elements of $D$.
By dovetailing we can define a total computable function
$\gamma:D^{<\omega}\to D$ such that, for all
$(d_0,...,d_k)\in D^{<\omega}$,
$$
d_i\not\leq\gamma(d_0,...,d_k)\mbox{ for all $i=0,...,k$}
$$
Let $b:\bbbn\to\words$ be such that $b(0)$ is the empty word
and $b(2n+1)=b(n)0$ and $b(2n+2)=b(n)1$. As is well known
(cf. Li \& Vitanyi \cite{livitanyi}, p.12),
$b$ is a total computable bijection which is length increasing:
$i<j\Rightarrow |b(i)|\leq|b(j)|$, so that
$$\{b_i:i\leq 2^k-2\}=\{q\in\words:|q|<k\}$$
Let $G=\max^\D g$ where $g:\words\times\bbbn\to D$ is total
computable.
Define a total computable $F:\WD$ as follows:
$$
F(p)=\gamma(g(b_0,0),...,g(b_{2^{2|p|}-2},0))
$$
By definition of $F$, we see that $g(q,0)\not\leq F(p)$ for all $q$
such that $|q|<2|p|$.
In particular, if $|q|<2|p|$ and $G(q)$ is defined,
since $g(q,0)\leq G(q)$ we have $F(p)\neq G(q)$.
This insures that $K_G(F(p))\geq 2|p|$.
Since, obviously, $K_F(F(p))\leq|p|$, we get
$K_G(F(p))\geq K_F(F(p))+|p|$.
Which proves that $K_G-K_F$ takes arbitrarily large values,
hence $G$ cannot be optimal in $\maxr[\WD]$.
Since $F$ is total computable, this also proves the last assertion
of Point 2.
\qed\end{pf}
Applying Theorem \ref{thm:KmaxRec} to $\bbbn$ and $\bbbz$
with the natural orderings, we get the following result.
It is interesting to compare Point 1 with
Proposition \ref{p:exitMach}.
\begin{cor}
1. The invariance theorem holds for the class
$\maxr[\words\to\bbbn]$.
Moreover, optimal functions in $\maxr[\words\to\bbbn]$ are optimal
for the class $\maxpr[\words\to\bbbn]$.
In particular, the Kolmogorov complexity associated to
$\maxr[\words\to\bbbn]$ coincides (up to a constant) with that
associated to $\maxpr[\words\to\bbbn]$.
\medskip\\
2. The invariance theorem fails for the classes
$\minr[\words\to\bbbn]$, $\maxr[\words\to\bbbz]$ and
$\minr[\words\to\bbbz]$.
\end{cor}
Since $\Reg$ with the inclusion ordering
(cf. \S\ref{sss:exQuotient}) has a minimum and a maximum element
(namely $\emptyset$ and $\widetilde{\Sigma}$), we get:
\begin{cor}
The invariance  theorem holds for the classes
$\maxr[\words\to\Reg]$ and $\minr[\words\to\Reg]$.
In particular, the associated Kolmogorov complexities
coincide (up to a constant) with those
associated to $\maxpr[\words\to\Reg]$ and $\minpr[\words\to\Reg]$.
\end{cor}

\end{document}